\documentclass[journal]{IEEEtran}
\usepackage{amsmath,amscd,amsbsy,amssymb,latexsym,url,bm,amsthm}
\allowdisplaybreaks
\usepackage{array}
\usepackage{algorithm}
\usepackage{algorithmic}
\usepackage{booktabs}
\usepackage{cite}
\usepackage{etoolbox} % to adjust some commands to highlight newly added references
\usepackage[T1]{fontenc}
\usepackage{graphicx}
\usepackage{multirow,makecell}
\usepackage{subeqnarray}
\usepackage[caption=false,font=footnotesize,subrefformat=parens,labelformat=parens]{subfig}
\usepackage{tabularx}
\usepackage{threeparttable}
\usepackage{tikz}
\usepackage{xcolor}
\usepackage{balance}
\usepackage{hyperref}
\hypersetup{hidelinks}

\newtheorem{assumption}{Assumption}
\newtheorem{theorem}{Theorem}
\newtheorem{lemma}[theorem]{Lemma}

\newtheorem{definition}{Definition}
\newtheorem{remark}{Remark}

\newcommand{\bigO}[1]{\ensuremath{\mathop{}\mathopen{}\mathcal{O}\mathopen{}\left(#1\right)}}
\newcommand{\overbar}[1]{\mkern 1.5mu\overline{\mkern-1.5mu#1\mkern-1.5mu}\mkern 1.5mu}
\newcommand{\ud}{\mathrm{d}}
\newcommand{\RNum}[1]{\uppercase\expandafter{\romannumeral #1\relax}}

\newcommand{\PR}[1]{\Pr\big\{#1\big\}}

\DeclareMathOperator*{\diam}{diam}

\newcolumntype{Y}{>{\centering\arraybackslash}X}

\makeatletter
\def\raisedotfill{%
  \leavevmode
  \cleaders \hb@xt@ .44em{\hss\raise0.5ex\hbox{.}\hss}\hfill
  \kern\z@}
\pretocmd\@bibitem{\color{black}\csname keycolor#1\endcsname}{}{\fail}
\newcommand\citecolor[1]{\@namedef{keycolor#1}{\color{blue}}}
\makeatother

\begin{document}
    \title{Private and Robust Distributed Nonconvex \\ Optimization via Polynomial Approximation}
    \author{Zhiyu He\textsuperscript{$\dagger$,$\ddagger$}~\IEEEmembership{Student Member, IEEE}, Jianping He\textsuperscript{$\dagger$}~\IEEEmembership{Senior Member, IEEE}, Cailian Chen\textsuperscript{$\dagger$}~\IEEEmembership{Member, IEEE}, \\ and Xinping Guan\textsuperscript{$\dagger$}~\IEEEmembership{Fellow, IEEE}
        \thanks{This work was supported in part by the National Natural Science Foundation of China under Grants 62373247, 92167205, 61933009, and 62025305 and in part by the Max Planck ETH Center for Learning Systems. This paper was presented in part at the 59\textsuperscript{th} IEEE Conference on Decision and Control\cite{he2020cdc}. \emph{(Corresponding author: Jianping He.)}

        % and the Max Planck ETH Center for Learning Systems
        \textsuperscript{$\dagger$}Department of Automation, Shanghai Jiao Tong University, Shanghai 200240, China, Key Laboratory of System Control and Information Processing, Ministry of Education of China, Shanghai 200240, China, and Shanghai Engineering Research Center of Intelligent Control and Management, Shanghai 200240, China (email: hzy970920@sjtu.edu.cn; jphe@sjtu.edu.cn; cailianchen@sjtu.edu.cn; xpguan@sjtu.edu.cn).

        \textsuperscript{$\ddagger$}Automatic Control Laboratory, ETH Z{\"u}rich, 8092 Z{\"u}rich, Switzerland (email: zhiyhe@ethz.ch).
        %Max Planck ETH Center for Learning Systems

        % \textsuperscript{\textsection} Max Planck Institute for Intelligent Systems, 72076 T{\"u}bingen, Germany.

        % Z.~He, J.~He, C.~Chen, and X.~Guan are with the Department of Automation, Shanghai Jiao Tong University, Shanghai 200240, China, also with the Key Laboratory of System Control and Information Processing, Ministry of Education of China, Shanghai 200240, China, and also with the Shanghai Engineering Research Center of Intelligent Control and Management, Shanghai 200240, China (email: hzy970920@sjtu.edu.cn; jphe@sjtu.edu.cn; cailianchen@sjtu.edu.cn; xpguan@sjtu.edu.cn).

        % Z.~He is also with the Automatic Control Laboratory, ETH Z{\"u}rich, 8092 Z{\"u}rich, Switzerland (email: zhiyhe@ethz.ch).
        
        % Republic of Korea, December 2020\cite{he2020cdc}.
        % This research work is partially sponsored by the National Key R\&D Program of China 2017YFE0114600, and NSFC of China (61973218, 61828301).
        }}

    % make the title area
    \maketitle

    \begin{abstract}
        There has been work that exploits polynomial approximation to solve distributed nonconvex optimization problems involving univariate objectives. This idea facilitates arbitrarily precise global optimization without requiring local evaluations of gradients at every iteration. Nonetheless, there remains a gap between existing guarantees and practical requirements, e.g., privacy preservation and robustness to network imperfections. To fill this gap and keep the above strengths, we propose a Private and Robust Chebyshev-Proxy-based distributed Optimization Algorithm (PR-CPOA). Specifically, to ensure both the accuracy of solutions and the privacy of local objectives, we design a new privacy-preserving mechanism. This mechanism leverages the randomness in blockwise insertions of perturbed vector states and hence provides a stronger privacy guarantee in the scope of ($\alpha,\beta$)-data-privacy. Furthermore, to gain robustness against network imperfections, we use the push-sum consensus protocol as a backbone and discuss its specific enhancements. Thanks to the purely consensus-type iterations, we avoid the privacy-accuracy trade-off and the bother of selecting proper step sizes in different settings. We rigorously analyze the accuracy, privacy, and complexity of the proposed algorithm. We show that the advantages brought by introducing polynomial approximation are maintained when all the above requirements exist.
    \end{abstract}


    % !TEX root = ..\article.tex
\section{Introduction}\label{sec:intro}
Distributed optimization enables multiple agents in a network to optimize the average of local objective functions in a collaborative manner. This global aim is achieved by exploiting local computations and communication between neighboring agents. Such a distributed architecture is preferable in various applications related to networked systems, e.g., distributed learning\cite{scutari2019distributed}, energy management\cite{yang2016distributed}, and resource allocation\cite{han2017differentially}. In these applications, the needs to improve efficiency, scalability, and robustness and protect privacy motivate the development of distributed strategies, which serve as plausible alternatives to their centralized counterparts. %\cite{nedic2018network}
%Distributed optimization enables multiple agents in a network to agree on the optimal points of the average of local objective functions.

\subsection{Motivations}
Considerable effort has been devoted to designing efficient primal\cite{nedic2009subgradient,shi2015extra} and dual-based\cite{makhdoumi2017convergence} distributed optimization algorithms and extending them to meet diverse practical requirements, including privacy preservation\cite{han2017differentially,hale2018cloud,nozari2018differentially,ding2022differentially,zhang2017dynamic}, time-varying directed communication\cite{nedic2017achieving,pu2020push}, and asynchronous computations to allow lack of coordination\cite{xu2017convergence,pu2020push}, delays, and packet drops\cite{tian2020achieving}. These extensions mainly focus on convex problems, and critical issues including privacy-accuracy trade-off\cite{nozari2018differentially} and network scaling\cite{nedic2017achieving} are explored.
%Most of these extensions focus on convex problems, and some critical issues including privacy-accuracy trade-off\cite{nozari2018differentially} and network scaling\cite{nedic2017achieving} are explored. 
% \rev{appropriate numbers of iterations}\cite{han2017differentially} and bounds for constant step-sizes\cite{pu2020push,tian2020achieving} are explored.
%shi2014linear,
%\cite{wu2017decentralized,tian2020achieving}

Despite their wide applicability, the above algorithms only ensure convergence to stationary points for nonconvex problems, and their loads of locally evaluating gradients or function values increase linearly with the number of iterations. These issues motivate the study of \cite{he2020distributed}, where polynomial approximations are introduced to substitute for univariate nonconvex local objectives, and a gradient-free, consensus-type iteration rule is adopted to exchange vectors of coefficients of local approximations. These designs help to achieve arbitrarily precise global optimization and reduce the costs of communication and local evaluations.
% More importantly, they separate the algorithm in \cite{he2020distributed} from typical gradient-based methods, and offer a new perspective to solve distributed optimization problems.
% Recently, \cite{he2020distributed} proposed a promising algorithm termed CPCA to address a class of constrained distributed nonconvex optimization problems. The core idea is to first use polynomial approximations (i.e., proxies) to substitute for general local objective functions, and then employ consensus protocols, where the information of coefficients of these proxies is exchanged, to enable agents to acquire a global proxy, and finally solve an easier approximate version of the original problem locally. The novel idea of employing polynomial approximation helps to achieve arbitrarily precise global optimization without demanding local evaluations of gradients or values of objective functions at every iteration.  

Nonetheless, there are two issues that limit the practical values of the algorithm CPCA (i.e., a Chebyshev-Proxy and Consensus-Based Algorithm) in \cite{he2020distributed}. First, it is not privacy-preserving due to the potential leakage of sensitive local objectives. This issue is critical because such leakage can cause the disclosure of secret local patterns. For instance, if local demand functions are revealed, then users' personal details (e.g., daily schedules) are at the risk of being inferred\cite{hale2018cloud}. %in smart grids, 
The above leakage stems from the consensus-type iterations, where vectors of coefficients of local approximations are directly exchanged. Once the adversaries obtain the exact initial vector of a target agent, they can recover a close estimate of its local objective. Hence, how to preserve the privacy of local objectives and quantify protection results is well worth consideration. Second, it only handles optimization over static undirected networks with perfect communication. Given that issues including time-varying directed links, lack of coordination, and packet drops are common in applications, it is meaningful to investigate their effects and find countermeasures to gain robustness. % against those network imperfections.
%The above issues lead to the study of this work.
In this paper, we aim to demonstrate that the introduction of polynomial approximation into distributed optimization allows enhancements to meet practical needs of privacy and robustness and, in the meantime, maintains advantages in solution accuracy and complexities.
%when privacy and robustness are taken into account.

\subsection{Contributions}
% In this paper, we exploit the idea of introducing polynomial approximation and 
We propose a Private and Robust Chebyshev-Proxy-based distributed Optimization Algorithm (PR-CPOA). % considering the requirements of privacy preservation and robustness to network imperfections. %, including time-varying directed communication and asynchrony. % due to lack of coordination, delays or packet drops.
The key idea is to construct Chebyshev polynomial approximations (i.e., proxies) for objectives, employ consensus-type iterations with a privacy-preserving and robust mechanism to exchange coefficient vectors of local proxies, and solve an approximate problem by optimizing the recovered global proxy. The main contributions are summarized as follows.

\begin{itemize}
    \item We propose PR-CPOA to solve distributed optimization problems with nonconvex univariate objectives and convex constraint sets, pursuing privacy preservation and robustness against network imperfections. We demonstrate that it maintains the advantages of CPCA in obtaining $\epsilon$ globally optimal solutions for any arbitrarily small given precision $\epsilon$ and being distributed terminable.
    \item We incorporate a new privacy-preserving mechanism to prevent sensitive local objective functions from being disclosed. This mechanism exploits two types of randomness. One lies in the obfuscation of local states with zero-sum random noises. The other is the randomness in the blockwise insertions of perturbed vector states to make their dimensions uncertain to the adversaries.
    %in both the obfuscation and the blockwise insertions of vector states. 
    %We thoroughly analyze the effect of privacy preservation through $(\alpha,\beta)$-data-privacy and demonstrate that the accuracy of solutions and the privacy of local objectives are simultaneously ensured.
    \item We prove that in the scope of $(\alpha,\beta)$-data-privacy\cite{he2018preserving}, a stronger privacy guarantee % (i.e., a lower disclosure probability) 
    is obtained compared to the design where existing algorithms (e.g., \cite{manitara2013privacy,mo2017privacy}) are directly extended to handle vector states. Moreover, we demonstrate that the solution accuracy and privacy of local objectives are simultaneously ensured. These guarantees are in contrast with differentially private distributed convex optimization algorithms\cite{han2017differentially,hale2018cloud,nozari2018differentially,ding2022differentially}.
    %he2019privacy
    % for the considered nonconvex problem
    \item We address the robustness issue in the face of various network imperfections. We employ the push-sum average consensus protocol\cite{kempe2003gossip} as a backbone to handle time-varying directed graphs and discuss its asynchronous extensions. We analyze the relationship between the accuracy of iterations and that of the obtained solutions, thus verifying that PR-CPOA remains effective against the above network imperfections. Thanks to the linear, consensus-type iterations, we are free from selecting appropriate step sizes in different settings. %which is a troublesome routine of gradient-based distributed methods.
    %We prove that the proposed algorithm keeps effective and accurate when such imperfections exist, and there is no need to carefully select proper step-sizes in different circumstances.
\end{itemize}

Compared to the conference version \cite{he2020cdc}, we \romannumeral1) explore the design and analysis of the privacy-preserving mechanism for approximation-based distributed optimization, \romannumeral2) investigate the strategies to deal with diverse network imperfections, and \romannumeral3) provide rigorous analysis of the accuracy and complexities.

\subsection{Organization}
%\emph{Organization:}
The remainder of this paper is organized as follows. Section~\ref{sec:formulation} describes the problem of interest and gives some preliminaries. Section~\ref{sec:algorithm} presents the algorithm PR-CPOA\@. Section~\ref{sec:analysis} analyzes the accuracy, privacy, and complexity of the proposed algorithm. Numerical evaluations are performed in Section~\ref{sec:experiment}, followed by the review of related work in Section~\ref{sec:literature}. Finally, Section~\ref{sec:conclusion} concludes this paper.
    % !TEX root = ..\article.tex
\section{Problem Description and Preliminaries}\label{sec:formulation}
Consider a network system of $N$ agents. Each agent $i$ owns a univariate local objective $f_i(x):X_{i} \to \mathbb{R}$ and a local constraint set $X_{i} \subset \mathbb{R}$. The network at time $t(t\in \mathbb{N})$ is described as a directed graph $\mathcal{G}^{t}=(\mathcal{V},\mathcal{E}^{t})$, where $\mathcal{V}$ is the set of agents, and $\mathcal{E}^{t} \subseteq \mathcal{V} \times \mathcal{V}$ is the set of edges. Note that $(i,j) \in \mathcal{E}^{t}$ if and only if (\textit{iff}) agent $j$ can receive messages from agent $i$ at time $t$. The script $k$ in parentheses denotes the index of components in a vector. We summarize important notations in Table~\ref{table:notations} for ease of reference.
%\rev{The superscript $t$, subscripts $i,j$, and script in parentheses $k$ denote the number of iterations, indexes of agents, and index of components in a vector, respectively.}

\begin{table}[!tb]
\renewcommand \arraystretch{.9}
% \small
% \setstretch{0.8}
\caption{Important Notations}
\label{table:notations}
\centering
    \begin{tabularx}{\columnwidth}{c X}
        \toprule
        \textbf{Symbol} & \textbf{Definition} \\
        \midrule
        % variables related to the topology
        %$\mathcal{G}$ & the network graph \\
        %$D$ & the diameter of $\mathcal{G}$ \\ % (i.e., the greatest distance between any pair of agents)
        %$U$ & an upper bound on $D$ known to all the agents \\
        % variables involved in the algorithmic design
        $c_{j}$ & the $j$-th Chebyshev coefficient \\
        %$T_{j}(u)$ & the $j$-th Chebyshev polynomial defined on $[-1,1]$ \\
        %$p_{i}^{0}$ & the initial local variable of agent $i$ \\
        %$p_{i}^{t}$ & the local variable of agent $i$ at time $t$ \\
        $f_i(x)$ & the local objective function of $f_i(x)$ \\
        $p_{i}^{(m_i)}(x)$ & a polynomial approximation of degree $m_i$ for agent $i$ \\
        $p_i^0$ & the coefficient vector of $p_{i}^{(m_i)}(x)$ for agent $i$ \\
        $m$ & the highest of the degrees of local approximations \\
        $\theta_i$ & a random noise vector used by agent $i$ \\
        $\Theta$ & the domain of a random variable $\theta_i(k)$ \\
        $g_{\theta_i(k)}(y)$ & the probability density function of a random variable $\theta_i(k)$ \\ 
        %$\mathbb{S}_{+}^{d}$ & the set of symmetric positive semidefinite matrices of order $d$ \\
        % variables related to the optimization problem
        %$\epsilon$ & the given solution accuracy (i.e., precision requirement) \\
        $\alpha$ & the estimation accuracy \\
        $\beta$ & the maximum disclosure probability \\
        $\mathcal{I}_{i}^{t}$ & the information set used by the adversaries to estimate agent $i$'s local objective at time $t$ \\ %$\mathcal{I}_{i}^{\textrm{own},t}$, $\mathcal{I}_{i}^{\textrm{in},t}$
        \bottomrule
    \end{tabularx}
\end{table}

\subsection{Problem Description}
We aim to solve the following constrained problem
\begin{equation}\label{problem:main_focus}
    \begin{split}
        \min_{x} \quad & f(x) = \frac{1}{N}\sum_{i=1}^{N} f_i(x) \\
        \textrm{s.t.} \quad & x \in X = \bigcap_{i=1}^{N} X_{i}
    \end{split}
\end{equation}
in a distributed, private, and robust manner. The global aim of optimization will be achieved through local communication and computations. Meanwhile, we will address practical requirements including preservation of the privacy of local objectives and robustness to time-varying directed communication and asynchrony. Some basic assumptions are as follows.

\begin{assumption}\label{assump:lipschitz_continous}
    The local objective $f_{i}(x)$ is Lipschitz continuous on $X_{i}$.
    %Every $f_{i}(x)$ is Lipschitz continuous on $X_{i}$.
\end{assumption}

\begin{assumption}\label{assump:constraint_set}
    The local constraint set $X_{i}$ is a closed, bounded, and convex set.
    %All $X_{i}$ are closed, bounded and convex sets.
\end{assumption}

Assumptions~\ref{assump:lipschitz_continous} and \ref{assump:constraint_set} are satisfied by problems of practical interests and are extensively made by the related literature (e.g., \cite{wai2017decentralized,scutari2019distributed} and the references therein).
%the literature on nonconvex distributed optimization (e.g., \cite{di2016next,wai2017decentralized,scutari2019distributed,chen2022distributed} and the references therein).
% Assumption~\ref{assump:lipschitz_continous} holds for a wide range of objective functions, including continuously differentiable functions. Both Assumptions~\ref{assump:lipschitz_continous} and \ref{assump:constraint_set} are commonly seen among the literature (e.g., \cite{scaman2017optimal,scutari2019distributed,xi2016distributed}) and are satisfied in problems of practical interests.

\begin{assumption}\label{assump:network_model}
     $\{\mathcal{G}^{t}\}$ is $B$-strongly-connected, i.e., there exists a positive integer $B$, such that for any $k\in \mathbb{N}$, the graph $\big(\mathcal{V}, \bigcup_{t=kB}^{(k+1)B-1} \mathcal{E}^{t}\big)$ is strongly connected.
\end{assumption}

Assumption~\ref{assump:network_model} states that the union graph is strongly connected for a time window of length $B$. It is weaker than that requiring connectivity at every time and is sufficient for information flow from one agent in networks to another\cite{nedic2017achieving}.

Problem~\eqref{problem:main_focus} involves nonconvex objective functions and convex constraint sets. % Therefore, it is a constrained nonconvex distributed optimization problem. 
Under Assumption~\ref{assump:constraint_set}, the set $X_{i}$ is a closed interval for any $i \in \mathcal{V}$. Hence, let $X_{i}=[a_{i},b_{i}]$, where $a_{i},b_{i} \in \mathbb{R}$. As a result, the intersection set $X$ is $[a,b]$, where $a=\max_{i \in \mathcal{V}} a_{i}$, $b=\min_{i \in \mathcal{V}} b_{i}$.

\subsection{Preliminaries}
% $\bullet$ 
\textit{Consensus Protocols:}~%\par
Let $\mathcal{N}_{i}^{\text{in},t} = \{j|(j,i)\in \mathcal{E}^{t}\}$ and $\mathcal{N}_{i}^{\text{out},t} = \{j|(i,j)\in \mathcal{E}^{t}\}$ be the sets of agent $i$'s in-neighbors %\footnote{As \cite{nedic2018network}, we assume that $i\in \mathcal{N}_{i}^{\text{in},t}, \forall t\in \mathbb{N}$, i.e., agent $i$ can always access its own information.}
and out-neighbors, respectively, and $d_{i}^{\text{out},t} = |\mathcal{N}_{i}^{\text{out},t}|$ (i.e., the cardinality of $\mathcal{N}_{i}^{\text{out},t}$) be its out-degree. Suppose that agent $i$ owns a local variable $x_{i}^{t}\in \mathbb{R}$. There are two consensus protocols, i.e., maximum consensus and average consensus, that allow agents to reach global agreement via local information exchange. The maximum consensus protocol is %\cite{saber2003consensus} is
\begin{equation}\label{eq:max_consensus}
    x_i^{t+1} = \max_{j \in \mathcal{N}_{i}^{\text{in},t}} x_j^{t}.
\end{equation}
It can be proven that with \eqref{eq:max_consensus}, all $x_i^{t}$ converge to $\max_{i\in \mathcal{V}} x_{i}^{0}$ in $T(\leq (N-1)B)$ iterations. %\cite{he2014time} %, i.e,
%\begin{equation*}
%    x_{i}^{t} = \max_{i\in \mathcal{V}} x_{i}^{0}, \quad \forall t \geq T,~i \in \mathcal{V}.
%\end{equation*}
The push-sum average consensus protocol\cite{kempe2003gossip} is
\begin{equation}\label{eq:push_sum}
    x_{i}^{t+1} = \sum_{j\in \mathcal{N}_{i}^{\text{in},t}} a_{ij}^{t} x_{j}^{t}, \quad y_{i}^{t+1} = \sum_{j\in \mathcal{N}_{i}^{\text{in},t}} a_{ij}^{t} y_{j}^{t},
\end{equation}
where $y_{i}^{t}\in \mathbb{R}$ is initialized to be $1$ for all $i\in \mathcal{V}$. The key to the convergence of \eqref{eq:push_sum} lies in constructing a column stochastic weight matrix $A^{t}\triangleq (a_{ij}^{t})_{N\times N}$. A feasible choice is to set $a_{ij}^{t}$ as $1/d_{j}^{\text{out},t}$ if $j\in \mathcal{N}_{i}^{\text{in},t}$, and as $0$ otherwise.
%of setting the weight $a_{ij}^{t}$ is %\cite{nedic2018network}
%\begin{equation}\label{eq:push_sum_weights}
%    a_{ij}^{t} =
%    \begin{cases}
%        1/d_{j}^{\text{out},t}, & \text{if } j\in \mathcal{N}_{i}^{\text{in},t}, \\
%        0, & \text{else.}
%    \end{cases}
%\end{equation}
In the implementation, every agent $j$ transmits $x_{j}^{t}/d_{j}^{\text{out},t}$ and $y_{j}^{t}/d_{j}^{\text{out},t}$ to its out-neighbors. With \eqref{eq:push_sum}, the ratio $z_{i}^{t} \triangleq x_{i}^{t}/y_{i}^{t}$ converges geometrically to the average of all the initial values $\overbar{x} = 1/N \sum_{i=1}^{N} x_{i}^{0}$\cite{kempe2003gossip}. %, i.e.,
%\begin{equation*}
%    \lim_{t\to \infty} z_{i}^{t} = \overbar{x}, \quad \forall i \in \mathcal{V}.
%\end{equation*}

% $\bullet$ 
\textit{Chebyshev Polynomial Approximation} %\par
%Chebyshev polynomial approximation 
focuses on using truncated Chebyshev series to approximate functions, thus facilitating numerical analysis. These series (i.e., approximations) are efficiently computed by interpolation. The degree $m$ Chebyshev interpolant $p^{(m)}(x)$ corresponding to a Lipschitz continuous function $g(x)$ defined on $[a,b]$ is 
\begin{equation}\label{eq:cheb_rep}
    p^{(m)}(x) = \sum_{j=0}^{m} c_{j} T_{j} \Big(\frac{2x-(a+b)}{b-a}\Big), \quad x\in [a,b],
\end{equation}
where $c_{j}$ is the Chebyshev coefficient, and $T_{j}(\cdot)$ is the $j$-th Chebyshev polynomial defined on $[-1,1]$ and satisfies $|T_{j}(x')| \leq 1, \forall x'\in [-1,1]$. As $m$ increases, $p^{(m)}(x)$ converges uniformly to $g(x)$ on the entire interval\cite{trefethen2013approximation}. %, i.e.,
%\begin{equation*}
%    \forall x \in [a,b],~\big|p^{(m)}(x)-g(x)\big| \rightarrow 0, \textrm{ as } m \rightarrow \infty.
%\end{equation*}
In practice, $p^{(m)}(x)$ with a moderate degree $m$ generally suffices to be a rather close approximation of $g(x)$\cite{trefethen2013approximation}. The dependence of $m$ on the smoothness of $g(x)$ and the specified precision $\epsilon$ is quantified in Sec.~\ref{subsec:complexity}. %Consequently, computing $p^{(m)}(x)$ becomes a practical way to construct an arbitrarily precise polynomial approximation for $g(x)$, as theoretically ensured by the \textit{Weierstrass Approximation Theorem}\cite[Theorem 6.1]{trefethen2013approximation}.

\subsection{Models of Adversaries of Privacy}
We consider \textit{honest-but-curious adversaries}\cite{wang2019privacy}. These adversaries are agents that faithfully follow the specified protocol but intend to infer information of the target agent $i$ based on the received data. They are linked with malicious or Byzantine adversaries\cite{su2021byzantine,kuwaranancharoen2020byzantine} sending manipulated information, in the sense that their existence will jeopardize the well-functioning of the network system. However, the difference is that here the goal is to prevent sensitive local information from being leaked instead of resilience to manipulation.
%This setup is different from malicious or Byzantine adversaries\cite{su2021byzantine,kuwaranancharoen2020byzantine} that perturb convergence by sending manipulated information. Moreover, the goal here is to prevent sensitive local information from being leaked instead of resilience to manipulation.

We are concerned with the issue of privacy disclosure arising in the consensus iterations of PR-CPOA\@. %In terms of these adversaries,
At time $t$, the exchanged information during iterations is
% For such iterations, at time $t$, the exchanged information serving as a basis for estimation consists of
\begin{equation*}
    \mathit{I}_{i}^{\textrm{own},t} = \{a_{ii}^{t},x_{i}^{t}\}, \quad \mathit{I}_{i}^{\textrm{in},t} = \{a_{ij}^{t},x_{j}^{t}|j\in \mathcal{N}_{i}^{\textrm{in},t}\},
\end{equation*}
i.e., the information sets of the states and weights of agent $i$ and those transmitted from $\mathcal{N}_{i}^{\textrm{in},t}$ to agent $i$, respectively. %As proven in \cite{mo2017privacy,he2019privacy}, the knowledge of $\bigcup_{t\in \mathbb{N}} \mathit{I}_{i}^{\textrm{own},t}$,~$\bigcup_{t\in \mathbb{N}} \mathit{I}_{i}^{\textrm{in},t}$ and the coupling between the locally added noises is a sufficient condition for the privacy compromise of noise-adding-based privacy-preserving consensus protocols. 
The following assumption characterizes the abilities of adversaries. %coupling relationship
\begin{assumption}\label{assump:bound_ability_adversary}
    At every time $t$, for the target agent $i$, the adversaries can always access $\mathit{I}_{i}^{\textrm{own},t}$ but can only obtain $\mathit{I}_{i}^{\textrm{in},t}$ with a probability whose upper bound is $p\in (0,1)$. % the full knowledge of
\end{assumption}
\begin{remark}
    We assume the constant access of $\mathit{I}_{i}^{\textrm{own},t}$ to include the scenario where some out-neighbors are adversaries and can therefore always receive the information transmitted by agent $i$, as in \cite{mo2017privacy,he2018preserving}. %The knowledge of $\mathit{I}_{i}^{\textrm{in},t}$ is assumed to be available with a probability not more than $p$ at time $t$. 
    For the assumption of $\mathit{I}_{i}^{\textrm{in},t}$, the rationality is that the switching nature of time-varying networks can inhibit the persistent and perfect access to $\mathit{I}_{i}^{\textrm{in},t}$. In practice, this setting holds if at time $t$, there exists a trustworthy agent whose link with agent $i$ occurs with a probability not less than $1-p$, or the adversaries are mobile and contact agent $i$ to gather $\mathit{I}_{i}^{\textrm{in},t}$ with a probability not more than $p$. %\cite{zhao2018resilient}.
    % out-neighbors $j \in \mathcal{N}_{i}^{\textrm{out},t}$
\end{remark}

\subsection{Privacy Requirement}\label{subsec:privacy_definition}
We consider preserving the privacy of agent $i$'s local objective $f_{i}(x)$, where $i=1,\ldots,N$. %As discussed in Sec.~\ref{sec:intro}, the knowledge of $f_i(x)$ may lead to the disclosure of sensitive local patterns, thus calling for protection.
In CPCA, local communication happens in its consensus iterations, where agents directly exchange and update their local variables $p_{i}^{0}$ that lie in a set $\mathcal{P} \subset \mathbb{R}^{m_i+1}$, i.e., $p_{i}^{0} \in \mathcal{P}$. These variables are the coefficient vectors of approximations $p_i^{(m_i)}(x)$ for $f_{i}(x)$. %%Once the adversaries obtain an estimation $\hat{p}_{i}$ of $p_{i}^{0}$, they will recover an approximation $\hat{f}_{i}(x)$ for $f_{i}(x)$, where $\hat{f}_{i}(x)$ is in the form of \eqref{eq:cheb_rep} with its coefficients in $\hat{p}_{i}$. 
Hence, $p_{i}^{0}$ is the sensitive information of $f_{i}(x)$ and its privacy should be preserved.
%% Hence, the way for adversaries to infer $f_{i}(x)$ is to obtain an estimation $\hat{p}_{i}^{0}$ for $p_{i}^{0}$ and then recover an approximation $\hat{f}_{i}(x):X\to \mathbb{R}$. The approximation $\hat{f}_{i}(x)$ is in the form of (\ref{eq:cheb_rep}) with its coefficients stored in $p_{i}^{0}$.
%We aim to design an average consensus algorithm that effectively preserves the privacy of $f_{i}(x)$, or more specifically, $p_i^0$. %This algorithm will be tailored to the case where agents own $p_i^0$ of different dimensions, given that the degrees of $p_i^{(m_i)}(x)$ vary.
To characterize the privacy effect, we use $(\alpha,\beta)$-data-privacy\cite{he2018preserving}. Let $\hat{p}_{i}$ be any estimation of $p_{i}^{0} \in \mathcal{P}$ based on the available information set $\mathcal{I}$. The definition of $(\alpha,\beta)$-data-privacy, where $\alpha \geq 0$ and $\beta \in [0,1)$, is given as follows.
\begin{definition}
    A distributed algorithm achieves ($\alpha,\beta$)-data-privacy for $p_{i}^{0} \in \mathcal{P}$ with a given $\mathcal{I}$ \textit{iff}
    \begin{equation}\label{eq:data_privacy_def}
        \max_{\hat{p}_i \in \mathcal{P}} \PR{\|\hat{p}_{i} - p_{i}^{0}\|_{1}\leq \alpha|\mathcal{I}} = \beta.
    \end{equation} 
\end{definition}
%\vspace{-1pc}
In \eqref{eq:data_privacy_def}, $\alpha$ and $\beta$ are parameters that indicate the estimation accuracy and the maximum disclosure probability of $p_{i}^{0}$, respectively. When $\alpha$ is specified, a smaller $\beta$ corresponds to a higher degree of privacy preservation. %The original definition of ($\alpha,\beta$)-data-privacy in \cite{he2018preserving} considers the estimation of scalar states, and here it is extended to handle vector states.

\begin{remark}
    (Connection with differential privacy (DP)\cite{han2017differentially}) DP emphasizes indistinguishability (i.e., a single substitution in the input will lead to a similar output). In some problems of network systems, however, we may further care about the exact privacy degree of the noise-adding mechanism in the face of estimation. Some adversaries may estimate the true local value, and the probability that such an estimation is close can be difficult to quantify through DP. Through the lens of data privacy\cite{he2018preserving}, we can explicitly characterize the relationship between estimation accuracy and disclosure probability.
\end{remark}
%\cite{he2020differential}

% \rev{Compared to other privacy metrics which emphasize indistinguishability (e.g., differential privacy\cite{he2020differential} and information-theoretic privacy\cite{gupta2020preserving}), data privacy offers a perspective of analyzing the relationship between estimation accuracy and maximum disclosure probability of sensitive information\cite{he2018preserving}.} 

%We use the $\ell_{1}$-norm of the error $\hat{p}_{i} - p_{i}^{0}$ to measure the estimation accuracy. This usage contributes to the neat relationship between the estimation accuracy of $p_{i}^{0}$ and that of $f_{i}(x)$, since $f_{i}(x)$ is closely approximated by $p_{i}^{(m_i)}(x)$, whose coefficients are stored in $p_{i}^{0}$. Detailed discussions are provided in Remark \ref{rem:privacy_vector}.
    % !TEX root = ..\article.tex
\section{Design of PR-CPOA}\label{sec:algorithm}
\begin{figure}[tb]
\begin{center}
    \includegraphics[width=0.9\columnwidth]{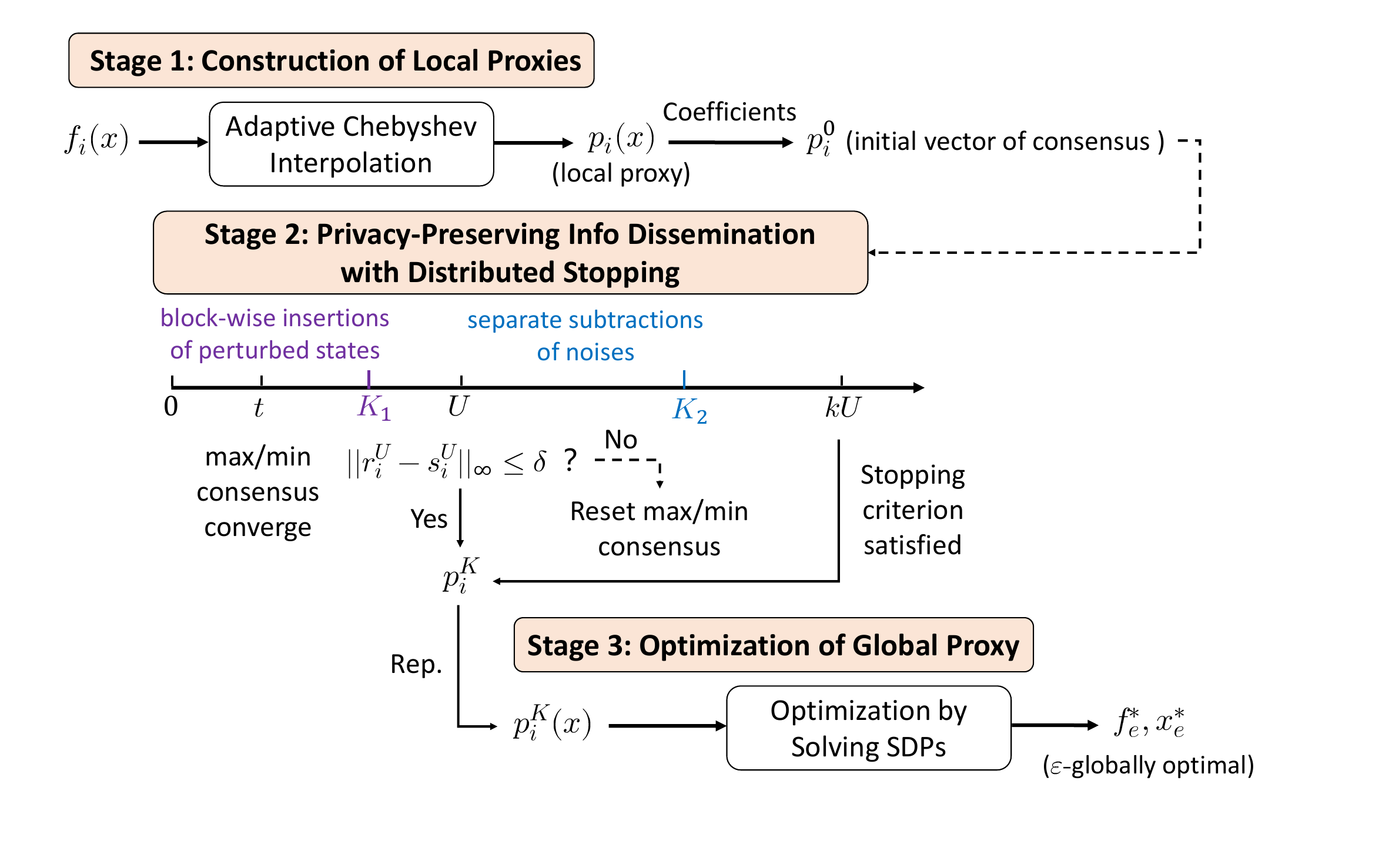}
    \caption{An overview of PR-CPOA.}
    \label{fig:overview}
    %\vspace{-4mm}
\end{center}
\end{figure}

We present the proposed PR-CPOA algorithm. It consists of three stages and is illustrated in Fig.~\ref{fig:overview}. %This algorithm consists of three stages and is illustrated in Fig.~\ref{fig:overview}.
% The proposed algorithm consists of three stages, whose details are discussed in the following three subsections.

% stage 1: initialization
\subsection{Construction of Local Chebyshev Approximations}\label{subsec:init}
In this stage, every agent $i$ computes a polynomial approximation $p_{i}^{(m_{i})}(x)$ of degree $m_{i}$ for $f_{i}(x)$ on $X=[a,b]$, s.t.
\begin{equation}\label{eq:approx_require}
   \big|f_{i}(x) - p_{i}^{(m_{i})}(x)\big| \leq \epsilon_{1}, \quad \forall x\in [a,b],
\end{equation}
where $\epsilon_{1}>0$ is a specified tolerance. This goal is achieved via the adaptive Chebyshev interpolation method\cite{boyd2014solving}. In this method, the degree of the interpolant is systematically increased until a certain stopping criterion is satisfied. First, agent $i$ sets $m_{i} = 2$ and evaluates $f_{i}(x)$ at the set $S_{m_{i}} \triangleq \{x_{0},\ldots,x_{m_{i}}\}$ of $m_{i}+1$ points by
\begin{equation}\label{eq:grid_point_evaluation}
    %\left\{
    %    \begin{aligned}
            x_{k} = \frac{b-a}{2}\cos\left(\frac{k\pi}{m_{i}} \right)+ \frac{a+b}{2}, \quad
            f_{k} = f_{i}(x_{k}),
    %    \end{aligned}
    %\right.
\end{equation}
where $k=0,1,\dots,m_{i}$. Then, it calculates the Chebyshev coefficients of the interpolant of degree $m_{i}$ by
\begin{equation}\label{eq:coeff_value_relation}
    c_{j} = \frac{1}{m_{i}} \left(f_{0} + f_{m_{i}}\cos(j\pi)\right) + \frac{2}{m_{i}} \sum_{k=1}^{m_{i}-1} f_{k}\cos\left(\frac{jk\pi}{m_{i}}\right),
\end{equation}
where $j=0,1,\dots,m_{i}$\cite{boyd2014solving}. At every iteration, the degree $m_{i}$ is doubled until the stopping criterion
\begin{equation}\label{eq:adapt_stop_rule}
    \max_{x_{k}\in \left(S_{2m_{i}} \setminus S_{m_{i}}\right)} \big|f_{i}(x_{k})-p_{i}^{(m_{i})}(x_{k})\big| \leq \epsilon_{1}
\end{equation}
is met, where $S_{2m_i} \setminus S_{m_i}$ is the set difference of $S_{2m_i}$ and $S_{m_i}$, and $p_{i}^{(m_i)}(x)$ takes the form of \eqref{eq:cheb_rep} with $\{c_{j}\}$ being its coefficients. Since $S_{m_i} \subset S_{2m_i}$, the evaluations of $f_{i}(x)$ are constantly reused. The intersection $X=[a,b]$ of local constraint sets is known by running some numbers of max/min consensus iterations as \eqref{eq:max_consensus} beforehand. 
% Then, the obtained $p_{i}^{m_{i}}(x)$ will satisfy the requirement \eqref{eq:approx_require} on accuracy\cite{boyd2014solving}.

% Remark on infinite or semi-infinite local constraint sets
%\begin{remark}
%    To handle problems where $X_{i}$ is an infinite or semi-infinite interval, we can use the technique of change-of-coordinate and choose rational Chebyshev functions as the basis for approximation, with technical assumptions on objectives\cite{boyd2014solving}. If local constraint sets are of mixed types (i.e., bounded, semi-infinite, or infinite), the best practice is to make agents run max/min consensus protocols in advance to know the intersection set, and then construct local approximations on a common basis. Additionally, when $a_{i}$ and $b_{i}$ cannot be obtained accurately (e.g., due to random noises), we can use the idea of sampling from scenario optimization\cite{garatti2019risk} and apply the proposed algorithm based on a random sample of constraints that are collected beforehand. The generalization property of the obtained solution can be characterized via the confidence bound for the constraint satisfaction\cite{garatti2019risk}.
%\end{remark}

% stage 2: iteration
\subsection{Privacy-Preserving Information Dissemination}\label{subsec:ppIterations}
Every agent now owns a local variable $p_{i}^{0}\in \mathcal{P} \subset \mathbb{R}^{m_i+1}$, which is the vector of coefficients of %local polynomial approximation
$p_{i}^{(m_i)}(x)$. In this stage, the goal is to enable agents to converge to the average $\overbar{p} = 1/N\sum_{i=1}^{N} p_{i}^{0}$ of their initial values\footnote{In this expression, those low dimensional vectors are extended with zeros when necessary to ensure the agreement in dimensions.} via a distributed mechanism, and the privacy of these initial values is preserved. %at the same time,
% Since $p_{i}^{0}$ is as an approximate representation of $f_{i}(x)$, the privacy preservation of $p_{i}^{0}$ implies that of $f_{i}(x)$.

We propose a privacy-preserving scheme to achieve this goal. The key ideas are \romannumeral1) adding random noises to $p_{i}^{0}$ to mask the true values, \romannumeral2) inserting the components of the perturbed initial states block by block, % to hide them within iterations,
thus making their dimensions uncertain to the adversaries, and \romannumeral3) subtracting the noises separately %in several randomly chosen rounds of iterations
to guarantee the convergence to the exact average. %The details are as follows.
% The backbone of this scheme is the push-sum average consensus protocol\cite{kempe2003gossip}.

$\bullet$ \textit{Step 1: Additions of Random Noises}\par
First, every agent $i$ generates a noise vector $\theta_{i}\in \Theta^{m_{i}+1}$, whose components are independent random variables in the domain $\Theta$. Then, it adds $\theta_{i}$ to its initial state $p_{i}^{0}$ to form a perturbed state $\tilde{p}_{i}^{0}$, i.e., $\tilde{p}_{i}^{0} = p_{i}^{0} + \theta_{i}$.
%\begin{equation*}
%    \tilde{p}_{i}^{0} = p_{i}^{0} + \theta_{i}.
%\end{equation*}

%\vspace{-0.9pc}
$\bullet$ \textit{Step 2: Blockwise Insertions of Perturbed States}\par
\begin{figure}[!tb]
    \centering
    \includegraphics[width=.5\columnwidth]{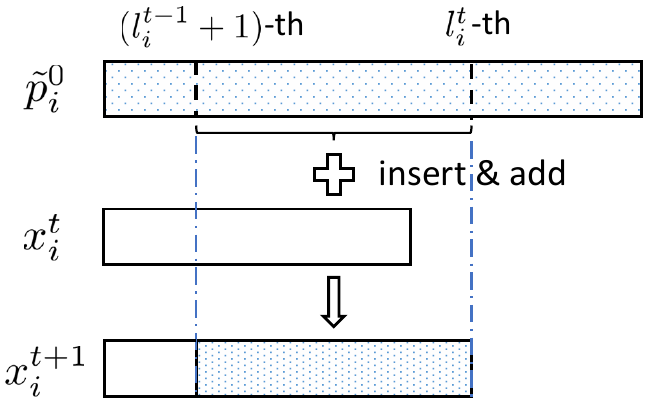}
    \caption{An illustration of blockwise insertions.}
    \label{fig:block_insert}
\end{figure}

Agents exchange and update their local variables $x_{i}^{t}$ and $y_{i}^{t}$ based on the push-sum consensus protocol\cite{kempe2003gossip}. The initial value of $y_{i}^{t}$ is all set as $1$ for all $i \in \mathcal{V}$. Instead of directly setting the initial value of $x_{i}^{t}$ as $\tilde{p}_{i}^{0}$, every agent $i$ will gradually extend $x_{i}^{t}$ with different blocks of $\tilde{p}_{i}^{0}$ in the first $K_{1}$ iterations. Let $(d_{i}^{1},\ldots,d_{i}^{K_{1}})$ be drawn from the multinomial distribution $\operatorname{multi}(m_{i}+1,1/K_1\cdot \mathbf{1}_{K_1})$, where $\mathbf{1}_{K_1}$ denotes an all-ones vector of size $K_1$. Then,
%\begin{equation*}
%    \sum_{t=1}^{K_{1}} d_{i}^{t} = m_{i}+1, \quad d_{i}^{t}\in \{0,\ldots,m_{i}+1\},~\forall t.
%\end{equation*}
$(d_{i}^{1},\ldots,d_{i}^{K_{1}})$ denote the numbers of components of $\tilde{p}_{i}^{0}$ that are inserted into $x_{i}^{t}$ at every iteration. Let $l_{i}^{0}=0$ and $l_{i}^{t}=\sum_{k=1}^{t} d_{i}^{k}$, where $t=1,\ldots,K_{1}$.
%\begin{equation*}
%    l_{i}^{0}=0, \qquad l_{i}^{t}=\sum_{k=1}^{t} d_{i}^{k},~~t=1,\ldots,K_{1}.
%\end{equation*}
At the $t$-th iteration, the $(l_{i}^{t-1}+1)$-th to $l_{i}^{t}$-th components of $x_{i}^{t}$ and $\tilde{p}_{i}^{0}$ are added together to form the corresponding components of $x_{i}^{t+}$. The remaining components of $x_{i}^{t}$ and $x_{i}^{t+}$ are the same. Specifically,
\begin{equation}\label{eq:insertion_rule}
    x_{i}^{t+}(k) =
    \begin{cases}
        x_{i}^{t}(k) + \tilde{p}_{i}^{0}(k),  & \text{for $k=l_{i}^{t-1}+1,\ldots,l_{i}^{t}$,} \\
        x_{i}^{t}(k), & \text{else}, \\
    \end{cases}
\end{equation}
where $x_{i}^{t}(k)$ and $\tilde{p}_{i}^{0}(k)$ denote the $k$-th components of $x_i^t$ and $\tilde{p}_i$, respectively. If the corresponding $x_{i}^{t}(k)$ is null (i.e., $k$ exceeds the size of $x_i^t$), then it is regarded as $0$. That is, in the first case of \eqref{eq:insertion_rule}, we add scalars and increase the sizes of vectors if necessary, thus avoiding disagreement in dimensions (see Fig.~\ref{fig:block_insert}). Then, agents transmit $x_{i}^{t+} \text{ and } y_{i}^{t}$ to their out-neighbors and update $x_{i}^{t+1} \text{ and } y_{i}^{t+1}$ by
\begin{equation}\label{eq:push_sum_insertion}
    x_{i}^{t+1}(k) = \sum_{j\in \mathcal{N}_{i}^{\text{in},t}} a_{ij}^{t} x_{j}^{t+}(k),\forall k, ~~ 
    y_{i}^{t+1} = \sum_{j\in \mathcal{N}_{i}^{\text{in},t}} a_{ij}^{t} y_{j}^{t}.
\end{equation}
The extension of $x_{i}^{t}$ is also involved in \eqref{eq:push_sum_insertion}. Hence, the dimension of $x_{i}^{t+1}$ is the same as the largest among $x_{j}^{t+}$, $j \in \mathcal{N}_{i}^{\text{in},t}$. At the end of the $K_{1}$-th iteration, all the components of $\tilde{p}_{i}^{0}$ have been inserted, and the size of $x_{i}^{t}$ is at least $m_{i}+1$.

$\bullet$ \textit{Step 3: Separate Subtractions of Noises}\par
In the following $K_{2}-K_{1}$ iterations, to ensure the convergence to the exact average $\overbar{p}$, every agent will properly subtract the added noises. Let $L_i \in \{1,\ldots,K_2 - K_1\} $ be a random integer such that
%between $1$ and $K_{2}-K_{1}$
\begin{equation}\label{eq:subtracted_noises}
    |\zeta_{i}(k)| > \alpha, \qquad \zeta_{i}(k) \triangleq \frac{\theta_{i}(k)}{L_i}.
\end{equation}
Note that $L_i$ can be drawn from various discrete distributions, e.g., the discrete uniform, binomial, and hypergeometric distributions. The choices of such distributions are up to the agents and are therefore unknown to the adversaries. For the $k$-th component of $x_{i}^{t}(k=1,\ldots,m_{i})$, at $L_i$ randomly selected numbers of iterations, every agent $i$ subtracts a fraction of the added noise $\zeta_{i}(k)$ from the updated state. That is, apart from performing \eqref{eq:push_sum_insertion}, it executes
\begin{equation}\label{eq:subtraction}
    x_{i}^{(t+1)+}(k) = x_{i}^{t+1}(k) - \zeta_{i}(k), ~~ \forall k = 1,\ldots, m_{i}.
\end{equation}
%\begin{subequations}\label{eq:subtraction}
%    \begin{align}
%        x_{i}^{t+1}(k) &= \sum_{j\in \mathcal{N}_{i}^{\text{in},t}} a_{ij}^{t} x_{j}^{t+}(k), \label{subeq:push_sum_x} \\
%        x_{i}^{(t+1)+}(k) &= x_{i}^{t+1}(k) - \zeta_{i}(k), ~~ \forall k = 1,\ldots, m_{i}, \\
%        y_{i}^{t+1} &= \sum_{j\in \mathcal{N}_{i}^{\text{in},t}} a_{ij}^{t} y_{j}^{t}.
%    \end{align}
%\end{subequations}
%The update in \eqref{subeq:push_sum_x} involves the extension of dimensions if necessary.
At the rest of the iterations, agents update their local variables by \eqref{eq:push_sum_insertion}, where $x_{i}^{t+}(k)$ is set as $x_{i}^{t}(k), \forall k, \forall t\geq K_{2}+1$.

$\bullet$ \textit{Auxiliary Step: Distributed Stopping}\par
To realize distributed stopping when the precision of iterations has met the requirement, after the $K_{2}$-th iteration, we utilize the max/min-consensus-based stopping mechanism in\cite{prakash2018distributed}. % The scheme in \cite{prakash2018distributed} deals with static digraphs, but it can be extended to handle time-varying digraphs, given that in this case the max/min consensus protocols still converge in finite time.
The following assumption is required by this mechanism.
\begin{assumption}\label{assump:upperbound}
    Every agent $i$ in $\mathcal{G} $ knows an upper bound $U$ on $(N-1)B$.
\end{assumption}
The bound $U$ can be obtained via the technique in \cite{martinez2018decentralized} to estimate $N$ and the prior knowledge of $B$. Specifically, there are two auxiliary variables, i.e., $r_{i}^{t}$ and $s_{i}^{t}$. They are initialized as $p_{i}^{K_{2}} = x_{i}^{K_{2}}/y_{i}^{K_{2}}$ and updated together with $x_{i}^{t}$ and $y_{i}^{t}$ by the following max/min consensus protocols
\begin{equation}\label{eq:max_consensus_auxiliary}
    r_{i}^{t+1}(k) = \max_{j\in \mathcal{N}_{i}^{\text{in},t}} r_{j}^{t}(k), ~~ s_{i}^{t+1}(k) = \min_{j\in \mathcal{N}_{i}^{\text{in},t}} s_{j}^{t}(k), ~~\forall k.
\end{equation}
With \eqref{eq:max_consensus_auxiliary}, $r_{i}^{t}(k)$ and $s_{i}^{t}(k)$ are guaranteed to converge in finite time to $\max_{j\in \mathcal{V}} p_j^{K_2}(k)$ and $\min_{j\in \mathcal{V}} p_j^{K_2}(k)$, respectively.
The required number of iterations for convergence is less than $(N-1)B$, and is therefore less than $U$. Hence, at time no later than $t=K_{2}+U$, all the local variables $x_{i}^{t}$, $y_{i}^{t}$, $r_{i}^{t}$, and $s_{i}^{t}$ become $(m+1)$-dimensional vectors, where $m \triangleq \max_{i\in \mathcal{V}} m_{i}$
%\cite{he2014time}
%\begin{equation*}
%    m \triangleq \max_{i\in \mathcal{V}} m_{i}
%\end{equation*}
is the maximum degree of all the local approximations. The variables $r_{i}^{t}$ and $s_{i}^{t}$ are reinitialized as $p_{i}^{t}$ every $U$ iteration to allow the continual dissemination of the recent information on $p_{i}^{t}$. When the stopping criterion
\begin{equation}\label{eq:itr_stop_rule}
    \|r_{i}^{K} - s_{i}^{K}\|_{\infty} \leq \delta \triangleq \frac{\epsilon_{2}}{m+1}
\end{equation}
is satisfied at the $K$-th iteration, agents terminate the iterations and set $p_{i}^{K} = x_{i}^{K}/y_{i}^{K}$.

% stage 3: optimization
\subsection{Polynomial Optimization by Solving SDPs}\label{subsec:PolyOpt}
In this stage, agents locally optimize the polynomial proxy $p_{i}^{K}(x)$ recovered from $p_{i}^{K}$ on $X$ to obtain $\epsilon$-optimal solutions of problem \eqref{problem:main_focus}. This problem is transformed into a semidefinite program (SDP) based on the sum-of-squares decomposition of non-negative polynomials\cite{blekherman2013semidefinite}. We refer the reader to Sec.~\RNum{3}-C of \cite{he2020distributed} for details on the transformed problem.

The transformed problems are SDPs and can therefore be efficiently solved via the primal-dual interior-point method\cite{boyd2004convex}. The iterations of this method are terminated when
% \begin{equation*}
    $0\leq f_{e}^{*} - p^{*} \leq \epsilon_{3}$,
% \end{equation*}
where $f_{e}^{*}$ is the obtained estimate of the optimal value $p^{*}$ of $p_{i}^{K}(x)$ on $X$, and $\epsilon_{3}>0$ is the specified precision. The optimal point $x_p^*$ of $p_{i}^{K}(x)$ on $X$ can then be calculated by the complementary slackness condition\cite{blekherman2013semidefinite}. %$X=[a,b]$
%The optimal points of $g_{i}^{K}(x)$ are computed from the complementary slackness condition\cite{blekherman2013semidefinite}. The optimal points of $p_{i}^{K}(x)$ on $X$ can then be calculated by \eqref{eq:opt_point_relation}.

% \subsection{Description of D-CPCA}
The full details of the proposed algorithm are summarized as Algorithm \ref{alg:alg_overall}. We set all the precision used in three stages, i.e., $\epsilon_{1},\epsilon_{2}\text{ and }\epsilon_{3}$, as $\epsilon/3$. Their sum equals $\epsilon$, thus helping to ensure the reach of $\epsilon$-optimality.
\begin{algorithm}[tb]
% \footnotesize
\small
\caption{PR-CPOA}
\label{alg:alg_overall}
% \setstretch{0.8}
    \begin{algorithmic}[1]
    \REQUIRE $f_{i}(x)$, $X_{i}=[a_{i},b_{i}]$, $U$, $K_1$, $K_2$, $\alpha$, $\epsilon$
    \ENSURE $f^{*}_{e}$, $x_p^*$. %for every agent $i \in \mathcal{V}$.
    \STATE {\bfseries Initialize:} $a_{i}^{0}=a_{i},b_{i}^{0}=b_{i},m_{i}=2$. \label{procedure:init}
    \FOR{{\bfseries each} agent $i\in \mathcal{V}$}
        \FOR{$t=0, \ldots, U-1$}
            \STATE $\displaystyle a_{i}^{t+1} = \max_{j\in \mathcal{N}_{i}^{\text{in},t}} a_{j}^{t}, ~ b_{i}^{t+1} = \min_{j\in \mathcal{N}_{i}^{\text{in},t}} b_{j}^{t}$.
        \ENDFOR
        \STATE Set $a=a_{i}^{t}$, $b=b_{i}^{t}$. \label{procedure:init_end}
        \\\raisedotfill % \textit{Construction of local proxies} \raisedotfill
        \STATE Calculate $\{x_{j}\}\textrm{ and }\{f_{j}\}$ by \eqref{eq:grid_point_evaluation}. \label{procedure:interpolation} \label{procedure:calc_lc_prxy}
        \STATE Calculate $\{c_{k}\}$ by \eqref{eq:coeff_value_relation}.
        \STATE If \eqref{eq:adapt_stop_rule} is satisfied (where $\epsilon_{1} = \frac{\epsilon}{3}$), go to step \ref{procedure:construct_vector}. Otherwise, set $m_{i} \leftarrow 2m_{i}$ and go to step \ref{procedure:interpolation}. \label{procedure:calc_lc_prxy_end}
        \\\raisedotfill % \textit{Push-sum consensus with distributed stopping}\raisedotfill
        \STATE Set $\tilde{p}_{i}^{0} = p_{i}^{0} + \theta_{i}$, $x_{i}^{0} = \text{null}$, $y_{i}^{0}=1$, $(d_{i}^{1},\ldots,d_{i}^{K_{1}})$ drawn from $\operatorname{multi}(m_{i}+1,1/K_1\cdot \mathbf{1}_{K_1})$, $l=1$, $\epsilon_2 = \frac{\epsilon}{3}$. \label{procedure:construct_vector}
        \FOR{$t=0,1,\ldots$}
            % the process of data insertion
            \IF{$t\leq K_{1}$}  %\COMMENT{Insertions of perturbed states}
                \STATE Extend $x_{i}^{t}$ to form $x_{i}^{t+}$ by \eqref{eq:insertion_rule}.
                \STATE Update $x_{i}^{t+1},\forall k \text{ and } y_{i}^{t+1}$ by \eqref{eq:push_sum_insertion}.
            \ELSIF{$K_{1}+1\leq t \leq K_{2}$}  %\COMMENT{Subtractions of noises}
                \FOR{{\bfseries each} component $k=1,\ldots,m_{i}$}
                    \STATE Update $x_{i}^{t+1}(k),\forall k \text{ and } y_{i}^{t+1}$ by \eqref{eq:push_sum_insertion}, or additionally by \eqref{eq:subtraction} if subtractions need to be performed.
                \ENDFOR
            \ELSE   %\COMMENT{Iterations with distributed stopping}
                \IF{$t=lU$}
                    \IF{$\|r_{i}^{t} - s_{i}^{t}\|_{\infty} \leq \epsilon_{2}/(m+1)$}
                        \STATE $p_{i}^{K} = x_{i}^{t}/y_{i}^{t}$. ~\textbf{break}
                    \ENDIF
                    \STATE $r_{i}^{t}=s_{i}^{t}=x_{i}^{t}/y_{i}^{t}$, $l \leftarrow l+1$.
                \ENDIF
                \STATE Update $x_{i}^{t+1}$, $y_{i}^{t+1}$, $r_i^{t+1}$, $s_i^{t+1}$ by \eqref{eq:push_sum} and \eqref{eq:max_consensus_auxiliary}.
            \ENDIF
        \ENDFOR
        % \STATE Set $p_{i}^{K} = x_{i}^{K}/y_{i}^{K}$. \label{procedure:average_consensus_end}
        \\\raisedotfill % \textit{Polynomial optimization by solving SDPs} \raisedotfill
        \STATE Solve the reformulated SDP with $\epsilon_{3}=\frac{\epsilon}{3}$ and return $f^{*}_{e}, x_p^*$. \label{procedure:poly_opt_slv}
    \ENDFOR
    \end{algorithmic}
\end{algorithm}

\section{Performance Analysis}\label{sec:analysis}
\subsection{Accuracy}
% We establish the accuracy of PR-CPOA. 
The following lemma guarantees the accuracy of the privacy-preserving iterations of the proposed algorithm.
\begin{lemma}\label{lem:itr_precision}
    If Assumptions 3 and 5 hold, when \eqref{eq:itr_stop_rule} is satisfied, we have 
    \begin{equation}\label{eq:itr_precision_requirement}
        \max_{i \in \mathcal{V}} \left \lVert p_{i}^{K} -\overbar{p} \right \rVert _{\infty} \leq \delta = \frac{\epsilon_{2}}{m+1}.
    \end{equation}
    %\rm \textbf{Proof.}
    %Please see Appendix~\ref{subsec:proof_itr_precision}. \qed
\end{lemma}

\begin{proof}
    % The proof relies on the investigation of consensus iterations and is referred to our online report\cite{he2021private}.
    Please see Appendix~\ref{subsec:proof_itr_precision}.
\end{proof}

% \begin{remark}
%     Lemma \ref{lem:itr_precision} is in the same form as \cite[Theorem 1]{he2020distributed}, but its proof is much more involved. Here we need to prove that with the insertions of perturbed data and separate subtractions of noises, the reach of exact average consensus is still ensured. We also need to verify the effectiveness of the stopping criterion (\ref{eq:itr_stop_rule}) in this case.
% \end{remark}

In the following theorem, we characterize the distance between the obtained solution $f_{e}^{*}$ and the optimal value $f^{*}$ of problem \eqref{problem:main_focus}, and the distance between the optimal point $x_{p}^{*}$ of $p_{i}^{K}(x)$ on $X$ (i.e., the returned solution) and the optimal point $x_{f}^{*}$ of problem \eqref{problem:main_focus}\footnote{Without loss of generality, we consider the case where $x_{f}^{*}$ is the single globally optimal point of problem~\eqref{problem:main_focus}. If there are multiple globally optimal points, we can perform a similar analysis by investigating the distance of $x_{p}^{*}$ to the set of all these globally optimal points.}.
% The following theorem demonstrates the accuracy of the proposed algorithm.
% We use $\epsilon \textrm{ and } f^{*}$ to denote the specified precision and the optimal value of problem (\ref{problem:main_focus}), respectively. 

\begin{theorem}\label{thm:alg_accuracy_result}
Suppose that Assumptions \ref{assump:lipschitz_continous}-\ref{assump:upperbound} hold. PR-CPOA ensures that every agent obtains $\epsilon$-optimal solutions $f_{e}^{*}$ for problem \eqref{problem:main_focus}, i.e., $\left|f_{e}^{*}-f^{*}\right| \leq \epsilon$.
%\begin{equation*}
%    \left|f_{e}^{*}-f^{*}\right| \leq \epsilon.
%\end{equation*}
Moreover,
\begin{equation*}
    |x_{p}^{*} - x_{f}^{*}| \leq \diam(S), \quad S = \left\{x\in X| f(x) \leq f(x_{f}^{*}) + {\textstyle \frac{4}{3}\epsilon}\right\}.
\end{equation*}
%\rm \textbf{Proof.}
%Please see Appendix~\ref{subsec:proof_accuracy}. \qed
\end{theorem}

\begin{proof}
    % The full proof is referred to our online report \cite{he2021private}. % due to space limit.
    Please see Appendix~\ref{subsec:proof_accuracy}.
\end{proof}

In Theorem~\ref{thm:alg_accuracy_result}, $\epsilon$ is any arbitrarily small specified precision, and $\diam(S)$ is the diameter of $S$, i.e., the maximum distance between any two points in $S$. It implies that for a given accuracy requirement, we can control the updates of the inner stages to ensure overall accuracy. The acquisition of $\epsilon$-optimal solutions benefits from the use of polynomial approximation and the characteristics of univariate objectives.

\subsection{Data-Privacy}
%Through $(\alpha,\beta)$-data-privacy\cite{he2018preserving}, 
We investigate the performance of PR-CPOA in preserving the privacy of $p_{i}^{0}$.
% we show that the developed algorithm preserves the privacy of $p_{i}^{0}$ and investigate the privacy-preserving property through the notion of data-privacy\cite{he2018preserving}.
We first define the information set $\mathcal{I}_{i}^{t}$ used by the adversaries at time $t$ for estimation. Let
\begin{align}
    \mathcal{I}_{i}^{t} &= \mathcal{I}_{i}^{\textrm{own},t}~ {\textstyle \bigcup} ~\mathcal{I}_{i}^{\textrm{in},t}, \label{eq:info_set} \\
    \mathcal{I}_{i}^{\textrm{own},t} &= {\textstyle \bigcup\limits_{s=1}^{t}} \mathit{I}_{i}^{\textrm{own},s} = {\textstyle \bigcup\limits_{s=1}^{t}} \{a_{ii}^{s},x_{i}^{s+}\}, \notag \\
    \mathcal{I}_{i}^{\textrm{in},t} &= {\textstyle \bigcup\limits_{s \in \mathbb{S}_{t}}} \mathit{I}_{i}^{\textrm{in},s} ~= {\textstyle \bigcup\limits_{s \in \mathbb{S}_{t}}} \{a_{ij}^{s},x_{j}^{s+}|j\in \mathcal{N}_{i}^{\textrm{in},s}\}. \notag
\end{align}
%\begin{equation*}
%    \mathcal{I}_{i}^{t} = \mathcal{I}_{i}^{\textrm{own},t}~ {\textstyle \bigcup} ~\mathcal{I}_{i}^{\textrm{in},t},
%\end{equation*}
%where
%\begin{alignat*}{2}
%    \mathcal{I}_{i}^{\textrm{own},t} &= {\textstyle \bigcup\limits_{s=1}^{t}} \mathit{I}_{i}^{\textrm{own},s} &&= {\textstyle \bigcup\limits_{s=1}^{t}} \{a_{ii}^{s},x_{i}^{s+}\}, \\
%    \mathcal{I}_{i}^{\textrm{in},t} &= {\textstyle \bigcup\limits_{s \in \mathbb{S}_{t}}} \mathit{I}_{i}^{\textrm{in},s} &&= {\textstyle \bigcup\limits_{s \in \mathbb{S}_{t}}} \{a_{ij}^{s},x_{j}^{s+}|j\in \mathcal{N}_{i}^{\textrm{in},s}\}.
%\end{alignat*}
The set $\mathbb{S}_{t}$ contains those numbers of iterations $s(s\leq t)$ when the adversaries have obtained the full knowledge of $\mathit{I}_{i}^{\textrm{in},s}$. Note that $\mathcal{I}_{i}^{t}$ consists of all the available information on the states and weights owned by and transmitted to agent $i$ up to the $t$-th iteration. Consider a random variable $X: \Omega \to \mathbb{R}$ whose distribution and any other relevant information are unknown. In this case, the probability that the adversaries can generate an accurate estimation $\hat{X}$ is small, i.e.,
\begin{equation}\label{eq:est_acc_no_info}
    \PR{|\hat{X}-X|\leq \alpha} \triangleq \gamma \ll p\max_{\nu \in \Theta} \int_{\nu-\alpha}^{\nu+\alpha} g_{\theta_{i}(k)}(y) \ud y, ~ \forall k,
\end{equation}
where $\Theta$ is the domain of the $k$-th component $\theta_i(k)$ of the noise vector $\theta_i$, and $g_{\theta_{i}(k)}(y)$ is the probability density function of $\theta_i(k)$\footnote{To illustrate, let $L_{\Omega}$ be the total length of $\Omega$. The sensible policy for the adversaries is to uniformly generate an estimation $\hat{X}$ from $\Omega$. Hence, $\gamma = 2\alpha/L_{\Omega}$. We can ensure \eqref{eq:est_acc_no_info} by choosing $g_{\theta_{i}(k)}(y)$ such that $\exists \tilde{\Theta} \triangleq [\tilde{v}-\alpha,\tilde{v}+\alpha], \forall y \in \tilde{\Theta}: g_{\theta_{i}(k)}(y) \gg 1/pL_{\Omega}$, which is not difficult to satisfy if $\Omega$ is a large domain (i.e., $L_{\Omega}$ is large).}.

Recall that we aim to preserve the privacy of $p_{i}^{0}\in \mathbb{R}^{m_{i}+1}$. Thanks to the blockwise insertions in \eqref{eq:insertion_rule}, the adversaries are unaware of the exact size $m_{i}+1$ of $p_i^0$. They do know the maximum size $m+1$, however, based on the received $p_{i}^{K}$. Hence, the estimation of $p_{i}^{0}$ consists of two parts, i.e., to estimate its components $p_{i}^{0}(k)$, where $k=1,\ldots,m_{i}+1$, and to infer that $p_{i}^{0}(k)$ is null for $k=m_{i}+2,\ldots,m+1$. Let $\alpha$ and $\alpha_{k}$ be the estimation accuracy of $p_{i}^{0}$ and each component $p_{i}^{0}(k)$, respectively, s.t.,
\begin{equation}\label{eq:element_accuracy_sum}
    \sum_{k=1}^{m_{i}+1} \alpha_{k} = \alpha, \quad \alpha_{k}\in [0,\alpha],~\forall k = 1,\ldots,m_{i}+1.
\end{equation}
Moreover, since $m_{i}$ is unknown and varies with $\epsilon_{1}$ and $f_{i}(x)$, it is viewed as a random variable by the adversaries. Let $F_{m_{i}|\mathcal{I}_{i}^{t}}(\cdot)$ be the cumulative distribution function of $m_{i}$ given $\mathcal{I}_{i}^{t}$. The following theorem characterizes the effects of privacy preservation of PR-CPOA\@.

\begin{theorem}\label{thm:privacy}
    If Assumptions \ref{assump:network_model} and \ref{assump:bound_ability_adversary} hold, given $\mathcal{I}_{i}^{t}$, PR-CPOA achieves $(\alpha,\beta)$-data-privacy for $p_{i}^{0}$, where $\{\alpha_{k}\}$ satisfies \eqref{eq:element_accuracy_sum},
    \begin{align}
        \beta &= \prod_{k=1}^{m_{i}+1} \beta_{k} \cdot \prod_{k=m_{i}+2}^{m+1} F_{m_i|\mathcal{I}_{i}^{t}}(k-2), \label{eq:data_privacy_beta} \\
        \beta_{k} &= \big(1-p^{K_{2}-K_{1}+1}\big) h_{i}(\alpha_{k}) + p^{K_{2}-K_{1}+1}, \label{eq:disclosure_prob_single_ele} \\
        h_{i}(\alpha_{k}) &= p\max_{\nu \in \Theta} \int_{\nu-\alpha_{k}}^{\nu+\alpha_{k}} g_{\theta_{i}(k)}(y) \ud y + \gamma. \notag
    \end{align}
    %\begin{equation}\label{eq:data_privacy_beta}
    %    \begin{split}
    %        \beta = \prod_{k=1}^{m_{i}+1} \beta_{k} \cdot \prod_{k=m_{i}+2}^{m+1} F_{m_i|\mathcal{I}_{i}^{t}}(k-2),
    %    \end{split}
    %\end{equation}
    %and
    %\begin{align*}
    %    \beta_{k} &= \prod_{k=1}^{m_{i}+1} \Big[\big(1-p^{K_{2}-K_{1}+1}\big) h_{i}(\alpha_{k}) + p^{K_{2}-K_{1}+1}\Big], \\
    %    h_{i}(\alpha_{k}) &= p\max_{\nu \in \Theta} \int_{\nu-\alpha_{k}}^{\nu+\alpha_{k}} f_{\theta_{i}(k)}(y) \ud y + \gamma.
    %\end{align*}
    %\rm \textbf{Proof.}
    %Please see Appendix~\ref{subsec:proof_privacy}. \qed
\end{theorem}

\begin{proof}
    Please see Appendix~\ref{subsec:proof_privacy}.
\end{proof}

Theorem~\ref{thm:privacy} states that PR-CPOA preserves the privacy of $p_{i}^{0}$ and characterizes such effects through ($\alpha,\beta$)-data-privacy. It clarifies the relationship between estimation accuracy and disclosure probability. We note that
\begin{equation*}
    \|\hat{p}_{i} - p_{i}^{0}\|_{1} = \sum_{k=1}^{m_{i}+1} |\hat{p}_{i}(k) - p_{i}^{0}(k)| \leq \sum_{k=1}^{m_{i}+1} \alpha_{k} = \alpha.
\end{equation*}
Hence, we sequentially consider the relationship between the estimation accuracy $\alpha_{k}$ and the maximum disclosure probability $\beta_{k}$ of each component $p_{i}^{0}(k)$, and then synthesize them to obtain the result concerning $p_{i}^{0}$. 
The interpretation of $\beta$ in \eqref{eq:data_privacy_beta} is as follows. Note that $\beta \in (0,1)$ is the product of a set of bounds $\beta_{k} \in (0,1)$ for disclosure probabilities of all the components $p_{i}^{0}(k)$ (see \eqref{eq:disclosure_prob_values}) and the probabilities (also in $(0,1)$) of correctly identifying null components (see \eqref{eq:disclosure_prob_null}). The bounds $\beta_{k}$ are derived via the law of total probability.
%Suppose that agent $i$ inserts its perturbed state $\tilde{p}_{i}^{0}(k)$ at time $s$. If the event that the adversaries know $I_{i}^{\text{in},t}$ for time $t=s-1$ and any time $t$ between $K_{1}+1$ and $K_{2}$ happens (the probability of which is not more than $p^{K_{2}-K_{1}+1}$), then the added noises $\theta_{i}(k)$ and the states $p_{i}^{0}(k)$ can be perfectly inferred. Otherwise, the disclosure probability will not exceed $h_{i}(\alpha_{k})$. The bounds $h_{i}(\alpha_{k})$ are derived likewise based on whether the adversaries know $I_{i}^{\text{in},s-1}$. If the adversaries know this information, then the maximum disclosure probability is
%% $\forall k=1,\ldots,m_{i}+1$
%\begin{equation*}
%    \max_{\nu \in \Theta} \int_{\nu-\alpha_{k}}^{\nu+\alpha_{k}} f_{\theta_{i}(k)}(y) \ud y,
%\end{equation*}
%which equals to the probability that the optimal distributed estimation falls into $[p_{i}^{0}(k)-\alpha_{k},p_{i}^{0}(k)+\alpha_{k}]$\cite{he2018preserving}. Otherwise, the disclosure probability is rather small since the adversaries own little relevant information of $p_{i}^{0}$.
The probabilities of the correct decision on null components are obtained based on whether the index $k$ exceeds $m_{i}+1$ (i.e., the dimension of $p_{i}^{0}$). %, which is viewed as a random variable by the adversaries).

If we directly extend existing algorithms \cite{manitara2013privacy,mo2017privacy,he2018preserving} to handle vector states $p_i^0$, then $\beta$ will at least equal $\prod_k \beta_k$. In contrast, the design of blockwise insertions causes adversaries to additionally identify null components, thus further reducing the disclosure probability of $p_{i}^{0}$. Additionally, such a benefit does not cause an increased communication complexity. %as we will see in Sec.~\ref{subsec:complexity}.
%he2019privacy

From \eqref{eq:data_privacy_beta}, we know that for $p_{i}^{0}$ of a larger size (i.e., with larger $m_{i}$), $\beta$ will generally be smaller, which implies a higher degree of privacy preservation. In addition, $\beta$ increases with $\alpha_{k}$ but decreases with $K_{2}-K_{1}$. These relationships support the intuitions that less accurate estimations can be acquired with higher probabilities, and more room for randomness leads to lower probabilities of privacy disclosure. %The choices of these parameters depend on the specific tolerances of estimation accuracy and numbers of iterations in applications.

\subsection{Further Discussions on Privacy and Robustness}\label{subsec:alg_dependability}
%We further discuss the performance and extension of the proposed algorithm, considering various requirements including privacy preservation and robustness to network imperfections. We compare PR-CPOA with other typical distributed optimization algorithms in Table~\ref{table:comparison_dependability}. The details are as follows.

$\bullet$ \textit{Privacy Guarantee.}
% We have shown in Theorem~\ref{thm:privacy} that PR-CPOA preserves the privacy of $p_{i}^{0}$ and analyzed the effects of preservation through ($\alpha,\beta$)-data-privacy. 
%Now, we study such effects via \textit{differential privacy}, 
%We now study the privacy effect via \textit{differential privacy}, which also uses a noise-adding mechanism and provides a strong privacy guarantee against adversaries owning arbitrarily much side information.
Although the characterization of data privacy in Theorem~\ref{thm:privacy} is different from differential privacy \cite{han2017differentially,nozari2017differentially} that emphasizes indistinguishability, they are closely related. The links lie in the use of noise-adding mechanisms to pursue privacy and the similar intuition that large noises ensure strong privacy preservation but degrade the utility of data. We now discuss the privacy effects via differential privacy.
%he2020differential
We define the database $D$ and the randomized query output $\mathcal{M}(D)$ until time $t$ as the set of initial states and the information set $\mathcal{I}_i^t$ of the adversaries (see \eqref{eq:info_set}), i.e.,
%\begin{equation*}
    $D = \{p_{i}^{0}|\forall i\in \mathcal{V}\}, ~ \mathcal{M}(D) = \mathcal{I}_i^t$,
%\end{equation*}
% \mathcal{M}(D) = \{x_{i}^{+}(t)|\forall t\in \mathbb{N},i\in \mathcal{V}\}
respectively. Based on \cite{nozari2017differentially}, in our setting, a privacy-preserving consensus protocol is \textit{$\epsilon$-differentially private} if
\begin{equation*}
    \PR{\mathcal{M}(D)\in \mathcal{O}}\leq e^{\epsilon}\PR{\mathcal{M}(D')\in \mathcal{O}}
\end{equation*}
holds for any $\mathcal{O}\subseteq \text{range}(\mathcal{M})$ and $\sigma$-adjacent $D,D'$ satisfying
\begin{equation*}
    \big\|p_{i}^{0}-(p_{i}^{0})'\big\|_{1}\leq
    \begin{cases}
        \sigma, & \text{if } i = i_{0}, \\
        0, & \text{if } i \neq i_{0}
    \end{cases}
\end{equation*}
for all $i\in \mathcal{V}$, where $i_{0}$ is some element in $\mathcal{V}$. Note that we have used correlated noises (see \eqref{eq:subtracted_noises}) to pursue the proximity of $p_{i}^{K}$ to the exact average $\overbar{p}$ (see Lemma~\ref{lem:itr_precision}), thus ensuring the accuracy of the obtained solutions (see Theorem~\ref{thm:privacy}). Based upon the impossibility result of simultaneously achieving exact average consensus and differential privacy\cite{nozari2017differentially}, we conclude that our algorithm is not $\epsilon$-differentially private. To pursue differential privacy at the cost of losing certain solution accuracy, we can add uncorrelated noises (e.g., independent Laplace noises) to the transmitted states at every iteration.

$\bullet$ \textit{Asynchrony.} %We discuss the asynchronous extension of the proposed algorithm.
We discuss the robustness issue from the angle of asynchronous paradigms.
%Compared to synchronous models, asynchronous paradigms 
%They are more desirable in applications for their increased efficiency in
They handle uncoordinated computations and imperfect communication, e.g., transmission delays and packet drops. The design of consensus-type information dissemination in Algorithm~\ref{alg:alg_overall} is synchronous. Its extension to cope with asynchrony is feasible and can benefit from the extensive research on asynchronous consensus protocols, including those allowing for random activations (e.g., gossip algorithms\cite{kempe2003gossip}), delays\cite{hadjicostis2013average}, packet drops, %\cite{bof2017average}
and all these issues\cite{tian2020achieving}. %In these protocols, the basic idea of proving convergence is to first transform asynchronous models to synchronous counterparts over augmented graphs, where virtual nodes and edges are added to facilitate the analysis, and then establish the convergence of synchronous models.
The aforementioned asynchronous protocols converge deterministically to the average of initial values. If they are incorporated into PR-CPOA, by Lemma~\ref{lem:itr_precision} and Theorem~\ref{thm:alg_accuracy_result}, the accuracy of the obtained solutions is still guaranteed, although the proof will be relatively more involved. In addition, since the iterations of Algorithm~\ref{alg:alg_overall} are consensus-based and gradient-free, there is no need to select varying step sizes in different circumstances of asynchrony.
% do not involve gradients

% $\bullet$ \textit{Complexity.}
\subsection{Complexity}\label{subsec:complexity}
We present a lemma about the dependence of the degree $m_{i}$ of the local approximation $p_{i}^{(m_{i})}(x)$ on the specified tolerance $\epsilon_{1}$ and the smoothness of the local objective $f_{i}(x)$.
\begin{lemma}[\hspace{1sp}\cite{he2020distributed}]\label{lem:degree_approx}
    If $f_{i}(x)$ and its derivatives through $f_{i}^{(v-1)}(x)$ are absolutely continuous and $f_{i}^{(v)}(x)$ is of bounded variation on $X_{i}$, then $m_{i} \sim \mathcal{O}(\epsilon_1^{-1/v})$. If $f_{i}(x)$ is analytic on $X_{i}$, then $m_{i} \sim \mathcal{O}(\ln\frac{1}{\epsilon_1})$.
\end{lemma}
Lemma~\ref{lem:degree_approx} suggests that for functions that are smooth to some extent, polynomial approximations of moderate degrees (e.g., of the order of $10^{1} \sim 10^{2}$) can serve as rather accurate representations\cite{trefethen2013approximation}.
The following theorem describes the complexities of PR-CPOA in terms of the maximum degree of local proxies $m$ and the solution accuracy $\epsilon$, which equal $\max_{i\in \mathcal{V}}$ and $3\epsilon_{1}$, respectively. We measure the computational complexity via the order of flops\footnote{A flop is defined as one addition, subtraction, multiplication or division of two floating-point numbers\cite{boyd2004convex}.} 
and use $F_{0}$ to denote the cost of flops in one evaluation of $f_{i}(x)$\footnote{This cost depends on $f_{i}(x)$ \cite{boyd2004convex} and hence is not explicitly specified.}.
\begin{theorem}\label{thm:complexity}
    PR-CPOA ensures that every agent obtains $\epsilon$-optimal solutions for problem \eqref{problem:main_focus} with $\bigO{m}$ evaluations of local objective functions, $\bigO{\log \frac{m}{\epsilon}}$ rounds of inter-agent communication, $\bigO{\sqrt{m}\log \frac{1}{\epsilon}}$ iterations of primal-dual interior-point methods, and $\bigO{m\cdot \max(m^{3.5}\log \frac{1}{\epsilon},F_{0})}$ flops.
\end{theorem}

%\rm \textbf{Proof.}
\begin{proof}
    Please see Appendix~\ref{subsec:proof_complexity}. %\qed
\end{proof}

We compare PR-CPOA with typical distributed optimization algorithms in Table~\ref{table:comparison_dependability}. These algorithms exhibit sublinear convergence for nonconvex problems, and the complexities of evaluations and communication are $\bigO{\frac{1}{\epsilon}}$. In comparison, PR-COPA is more efficient in the complexities of function evaluations and communication rounds.

\begin{table*}[!tb]
%\tiny
% \scriptsize
% \small
% \setstretch{.8}
\centering
\renewcommand \arraystretch{0.6}
\begin{threeparttable}
\caption{Comparisons of PR-CPOA and Other Typical Distributed Optimization Algorithms}
\label{table:comparison_dependability}
    \begin{tabularx}{\linewidth}{c Y *{2}{c} *{3}{Y} c}
        \toprule
        {\bfseries \multirowcell{2}{Algorithms}} & {\bfseries \multirowcell{2}{Nonconvex \\ Objectives}} & \multicolumn{2}{c}{\bfseries Networks} & {\bfseries \multirowcell{2}{Privacy \\ Guarantee}} & {\bfseries \multirowcell{2}{Asynchrony}} & {\bfseries \multirowcell{2}{Accuracy \\ Guarantee}} & {\bfseries \multirowcell{2}{Complexities}} \\
        \cmidrule{3-4}
        & & Time-varying & Digraph & & & & \\
        %\midrule
        %Push-DIGing\cite{nedic2017achieving} & & \checkmark & \checkmark & & & \checkmark & scvx\tnote{1}: linear \\
        %\midrule
        %G-Push-Pull\cite{pu2020push} & & & \checkmark & & \checkmark & mean-square & scvx: linear \\
        %\midrule
        %SONATA\cite{scutari2019distributed} & \checkmark & \checkmark & \checkmark & & & \checkmark & \makecell{scvx: linear \\ ncvx\tnote{2}: $\bigO{\frac{1}{\epsilon}}$\tnote{3}} \\
        \midrule
        ASY-SONATA\cite{tian2020achieving} & \checkmark & & \checkmark & & \checkmark & \checkmark & \makecell{scvx\tnote{1}: linear \\ ncvx\tnote{1}: $\bigO{\frac{1}{\epsilon}}$\tnote{2}} \\
        %\midrule
        %Algorithm in \cite{han2017differentially} & & \multicolumn{2}{c}{Cloud-based} & DP\tnote{4} & & trade-off\tnote{5} &  \\
        %\midrule
        %Algorithm in \cite{hale2018cloud} & & \multicolumn{2}{c}{Cloud-based} & DP & & \dittostraight \tnote{6} &  \\
        \midrule
        Algorithm in \cite{nozari2018differentially}\tnote{3} & & \checkmark & \checkmark & DP & \checkmark & trade-off\tnote{4} & \\  %\makecell{\revise{square} \\ \revise{integrable}}
        \midrule
        %\makecell{Subspace \\ perturbation\cite{li2020privacy}} & & & \checkmark & \makecell{mutual- \\ information- \\ privacy} & & \checkmark & \\
        %\midrule
        FS protocol\cite{gupta2020preserving} & & & & \makecell{statistical \\ privacy} & & \checkmark & scvx: $\bigO{\frac{1}{\epsilon}}$ \\
        \midrule
        {\bfseries PR-CPOA} & \checkmark & \checkmark & \checkmark & \makecell{($\alpha,\beta$)- \\ data-privacy \\ (Theorem \ref{thm:privacy})} & \checkmark\tnote{5} & \checkmark & 
            \makecell{$0^{\textrm{th}}$-ord.\ oracle: $\bigO{m}$ \\ Commn.: $\bigO{\log \frac{m}{\epsilon}}$ \\ (Theorem~\ref{thm:complexity})} \\ %PD itr.: $\bigO{\sqrt{m}\log \frac{1}{\epsilon}}$} \\
        \bottomrule
        \addlinespace[0.5ex]
    \end{tabularx}
    \begin{tablenotes}[para]\scriptsize
        \item[1] ``scvx'' and ``ncvx'' refer to strongly-convex and nonconvex objective functions, respectively.
        \item[2] The complexities of inter-agent communication and evaluations of gradients (i.e., queries of the first-order oracle) are $\bigO{1/\epsilon}$.
        %\item[2] The convergence time is $\bigO{1/\epsilon}$, implying that the complexities of inter-agent communication and evaluations of gradients (i.e., queries of the first-order oracle) are $\bigO{1/\epsilon}$.
        \item[3] This strategy of function perturbation can be combined with any distributed convex constrained optimization algorithms. We place $\checkmark$ to imply feasibility.
        %to ensure differential privacy
        %The authors of \cite{nozari2018differentially} proposed a general strategy of function perturbation to ensure differential privacy. This strategy can be combined with any distributed convex constrained optimization algorithms to take effect. Hence, we place $\checkmark$ to some blocks in this row to imply feasibility.
        %\item[4] ``DP'' stands for differential privacy.
        \item[4] There is a trade-off between accuracy and privacy.
        \item[5] Detailed discussions are provided in Sec.~\ref{subsec:alg_dependability}.
        %\item[7] See Lemma~\ref{lem:degree_approx} and Theorem~\ref{thm:complexity} for details.
    \end{tablenotes}
\end{threeparttable}
\end{table*}

\begin{remark}
    Although PR-CPOA involves the exchange of $m$-dimensional vectors, its total transmission costs in communication can be acceptable given i) the decreased rounds of communication (see Theorem~\ref{thm:complexity} and Table~\ref{table:comparison_dependability}) and ii) the typically moderate degrees $m$ of approximations in numerical practice (see \cite{trefethen2013approximation} and also the discussion below Lemma~\ref{lem:degree_approx}).
\end{remark}

%\subsection{Discussions on Multivariate Extensions}
\subsection{Applications and Extensions}
In this paper, we consider problems with univariate objectives to highlight the advantages brought by exploiting polynomial approximation, e.g., achieving efficient optimization of nonconvex problems and allowing for enhancement to be private and robust when diverse practical needs exist. We discuss the relevant applications and multivariate extensions.

$\bullet$ \textit{Application Scenario.}
Some multivariate distributed optimization problems are naturally separable or can be transformed into a separable form through a change of coordinates. Consider the regularized facility location problem\cite{boyd2004convex,metel2019simple}
\begin{equation}\label{eq:facility_location}
    \min_{x \in X^n} ~ \sum_{i=1}^{N} w_i \|x-u_i\|_1 + \tilde{g}(x).
\end{equation}
In problem~\eqref{eq:facility_location}, $X^n \subset \mathbb{R}^n$ is a constraint set, $u_i \in \mathbb{R}^n$ is the location of agent $i$, $w_i$ is a nonnegative weight, and $\tilde{g}(x) \triangleq \sum_{i=1}^{n} \kappa \log(1+|x(k)|/\nu)$ is the log-sum penalty, where $\kappa,\nu>0$ are parameters. Note that $\tilde{g}(x)$ is one of the common surrogates of the $\ell_0$ norm to promote sparsity, and it is a nonconvex Lipschitz continuous regularizer\cite{metel2019simple}. Problem~\eqref{eq:facility_location} can be solved by considering a set of univariate subproblems as follows
%\vspace{-1em}
\begin{equation*}
    \min_{x(k) \in X} ~ \sum_{i=1}^{N} \left(w_i |x(k) - u_i(k)| + \kappa \log(1+|x(k)|/\nu) \right),
\end{equation*}
where $k = 1,\ldots,n$. The proposed algorithm applies to the above problem involving nonconvex Lipschitz objectives.

The key idea of exploiting function approximation is also of interest to address other problems of data analytics in network systems. We can use function approximation to represent local features and perform consensus to obtain the global feature (e.g., the average, product, or other functions of local features), thus facilitating operations of statistics or estimation\cite{he2020distributed}.
%li2021consensus

$\bullet$ \textit{Multivariate Extension.}
%We now discuss the multivariate extension of the proposed algorithm. 
Compared to the univariate setup, the differences of the multivariate extension mainly lie in the stages of initialization and optimization of approximations. Specifically, let $L_{2}(X)$ be the set of square-integrable functions over $X\subset \mathbb{R}^{n}$ and $f_{i}(x)\in L_{2}(X)$ be a general local objective. Then, there exists an orthonormal basis $\{h_{k}(x)\}_{k\in \mathbb{N}_{+}}$ (e.g., orthonormalization of Taylor polynomials) and an arbitrarily precise approximation
%\begin{equation*}
    $\hat{f}_{i}(x) = \sum_{k=1}^{m} c_{k} h_{k}(x)$
%\end{equation*}
for $f_{i}(x)$, where $\{c_{k}\}_{k=1}^{m}$ is the set of coefficients. Afterward, agents exchange and update their local vectors that store these coefficients (as in Sec.~\ref{subsec:ppIterations}) and acquire an approximation for the global objective. Finally, they locally optimize this approximation via the tools for polynomial optimization or for finding stationary points of general nonconvex functions\cite{zhang2020complexity}. % thus obtaining desired solutions. 
The main challenge of this extension lies in systematically calculating $\hat{f}_{i}(x)$ of an appropriate degree to strictly satisfy the accuracy requirement. It calls for further investigation and careful analysis and is still among our ongoing work.
% The aforementioned idea of extensions calls for further investigation and careful analysis and is still among our ongoing work.
Nonetheless, the idea of introducing approximation still i) offers a new perspective of representing possibly complicated objectives by simple coefficient vectors of approximations, thus allowing the dissemination of the global objective in the network, ii) naturally allows using random blockwise insertions to pursue a strong privacy guarantee in data privacy, and iii) is promising to obtain approximate globally optimal solutions in a distributed and asymptotic manner.

%\subsection{Robustness}
%In this section, we discuss the robustness of the developed algorithm in face of various issues arising in applications. These issues mainly include packet drops, transmission delays and asynchronous computations.

%In terms of D-CPCA, the aforementioned issues will influence the well-functioning of its inner consensus iterations, thus degrading the performance of the whole algorithm. Fortunately, there has already been extensive research on designing robust consensus algorithms that properly handle these issues. Hence, these designs can be incorporated into D-CPCA to further improve its robustness. Specifically, random packet drops can be equivalently viewed as random broken links that result in time-varying topologies, and thus D-CPCA can still remain effective. Some delays that exceed a certain threshold can be regarded as packet drops, while others that are smaller can be handled by the elegant technique of fusing possibly delayed neighboring states as soon as they arrive\cite{hadjicostis2013average}. As for realizing asynchronous computations, we can turn to gossip-type average consensus algorithms, where only a single link between neighboring agents is activated at every time slot and no synchronization of steps is needed\cite{boyd2006randomized}.

%may also discuss the influence of the insertions and subractions on the convergence rates of consensus iterations.
    % !TEX root = ..\article.tex
\section{Numerical Evaluations}\label{sec:experiment}
% We perform numerical experiments to illustrate the performance of PR-CPOA\@.
Consider a network with $N=20$ agents. At each time $t$, besides itself, every agent $i$ has two out-neighbors. One belongs to a fixed cycle, and the other is chosen uniformly at random. Hence, $\{\mathcal{G}^{t}\}$ is $1$-strongly-connected. We set the local constraint sets as the same interval $X=[-1,1]$ and generate the local objective function $f_{i}(x)$ of agent $i$ by
\begin{equation*}%\label{eq:local_obj}
    f_{i}(x) = \frac{a_{i}}{1+e^{-x}} + b_{i}\log (1+x^{2}),
\end{equation*}
where $a_{i}\sim \mathcal{N}(10,2) \text{ and } b_{i}\sim \mathcal{N}(5,1)$ are independent and identically distributed samples drawn from normal distributions. It follows that $f_{i}(x)$ is nonconvex and Lipschitz continuous on $X$. We use the Chebfun toolbox\cite{trefethen2013approximation} to construct Chebyshev polynomial approximations $p_{i}(x)$ for the local objective functions $f_{i}(x)$.
% Figure \ref{fig:performance-D-CPOA} presents the performance of D-CPOA\@.

\begin{figure*}
    %\centering
    %\subcaptionbox{Convergence\label{fig:convergence-PR-CPOA}}
    %{\includegraphics[width=0.3\linewidth]{figure/precision.eps}} \hfil
    %\subcaptionbox{($\alpha,\beta$)-data-privacy\label{fig:dataPrivacy}}
    %{\includegraphics[width=0.3\linewidth]{figure/dataPrivacy.eps}} \hfil
    %\subcaptionbox{Robustness\label{fig:robustness}}
    %{\includegraphics[width=0.3\linewidth]{figure/robustness.eps}}
    %\caption{Performance of PR-CPOA\@.}
    %\label{fig:performance-PR-CPOA}

    \centering
    \subfloat[]{\includegraphics[width=0.3\linewidth]{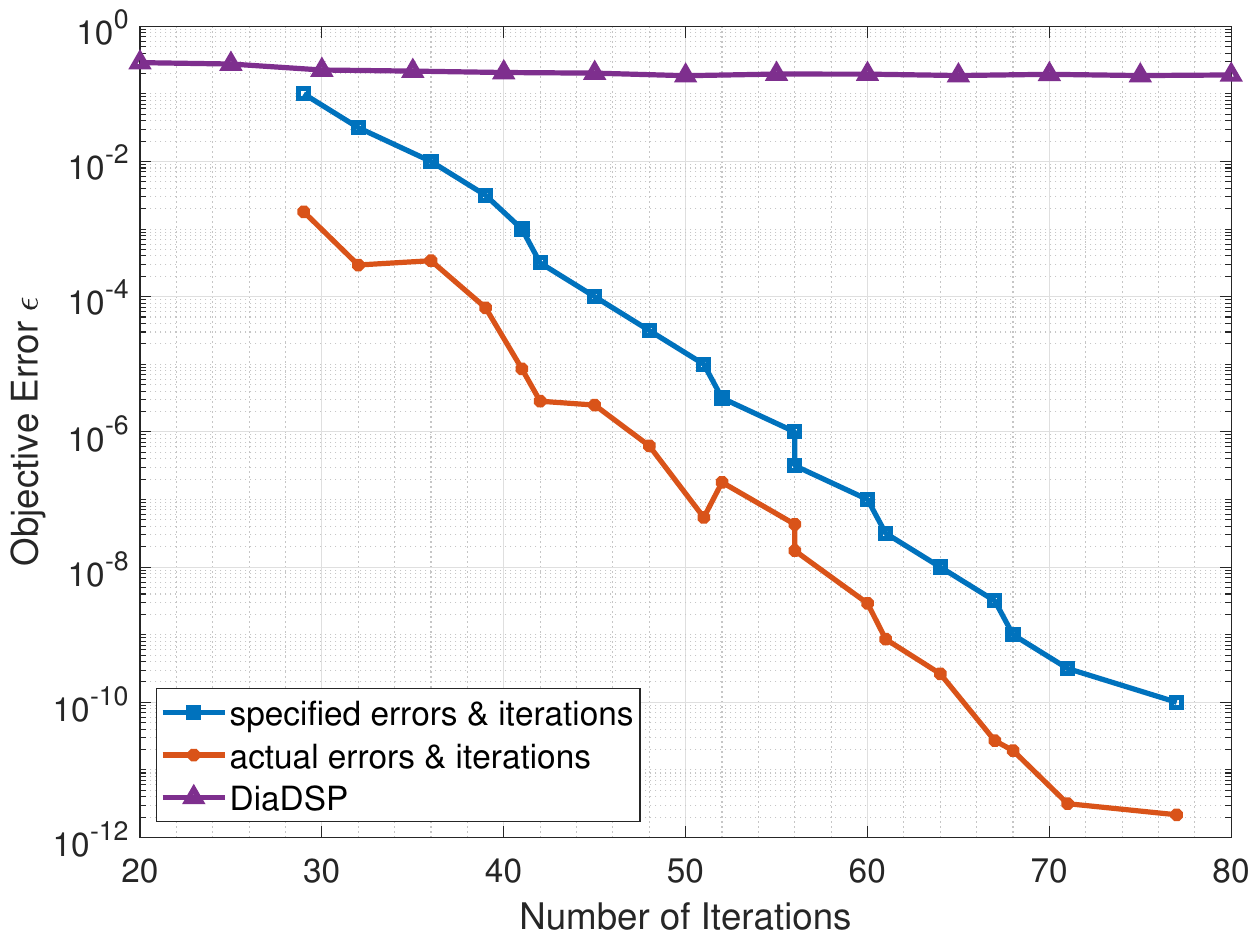}\label{fig:convergence-PR-CPOA}} \hfil
    \subfloat[]{\includegraphics[width=0.3\linewidth]{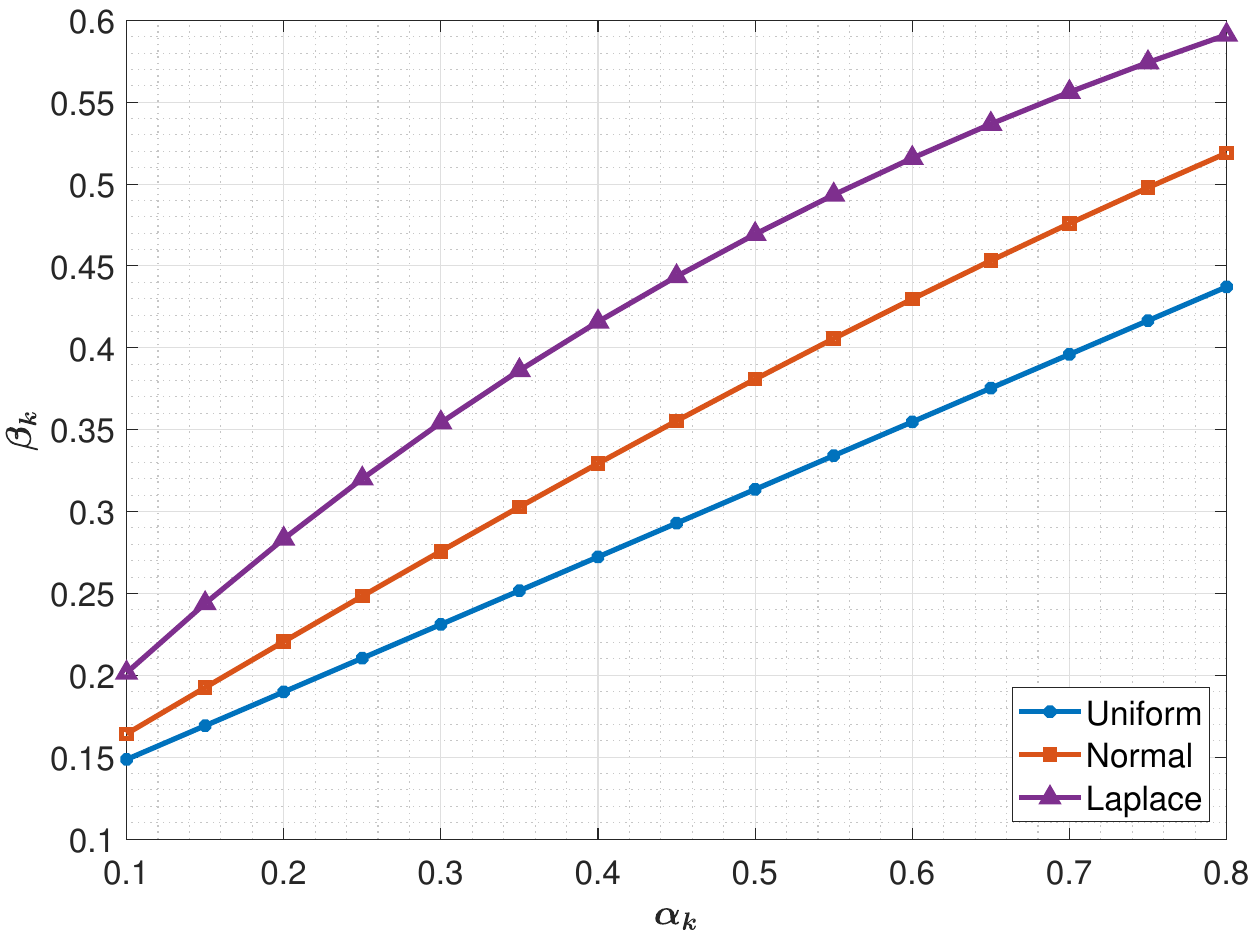}\label{fig:dataPrivacy}} \hfil
    %\subfloat[]{\includegraphics[width=0.25\linewidth]{figure/ppConsensusError.eps}\label{fig:ppConErr}} \hfil
    \subfloat[]{\includegraphics[width=0.3\linewidth]{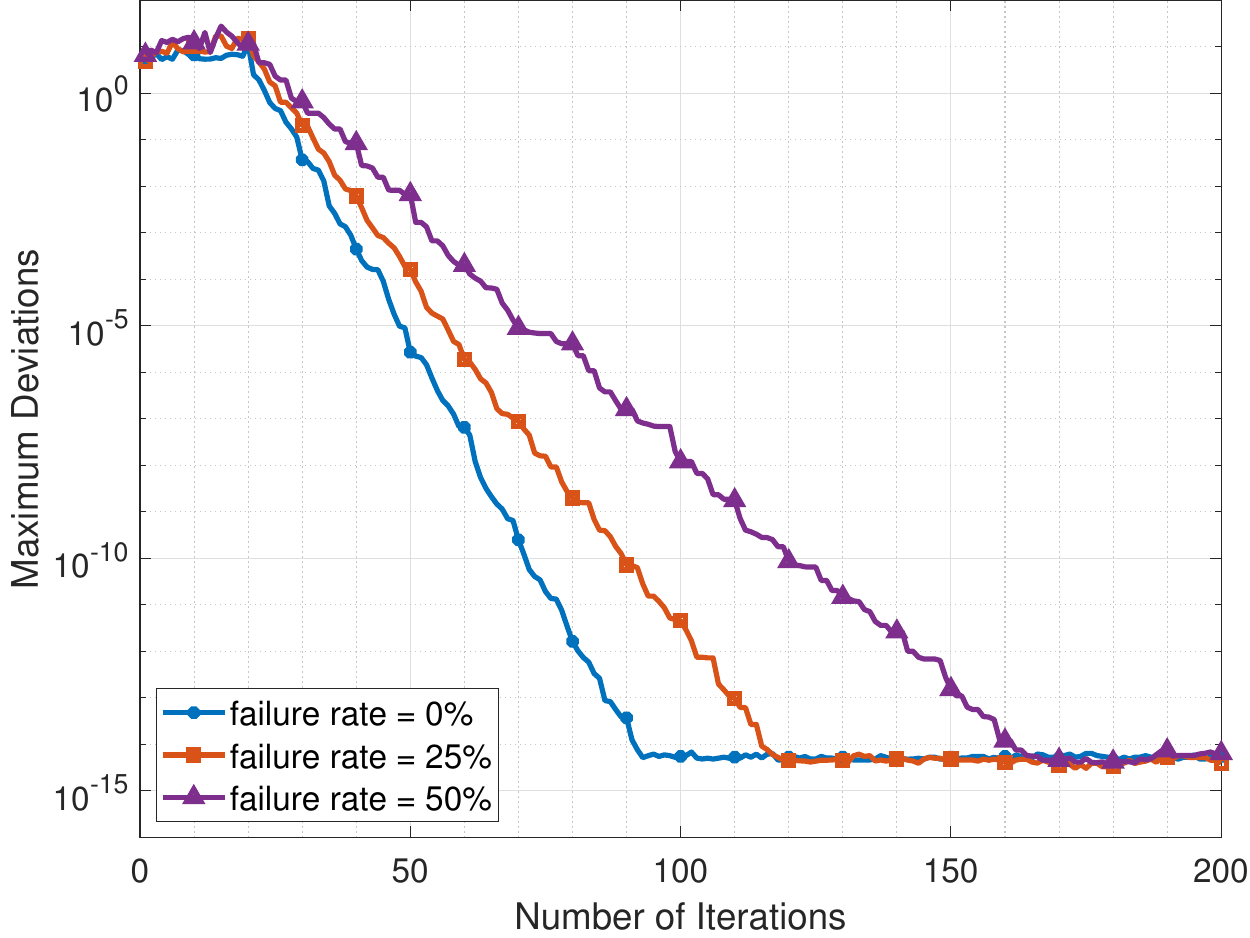}\label{fig:robustness}}
    \caption{Performance of PR-CPOA\@. (a) Convergence. (b) ($\alpha,\beta$)-data-privacy. %(c) Convergence of the proposed protocol in Sec.~\ref{subsec:ppIterations} and other privacy-preserving consensus protocols. 
    (c) Robustness.}
    \label{fig:performance-PR-CPOA} 
     \vspace{-1ex}
\end{figure*}

For PR-CPOA, we set $K_{1}=10,K_{2}=20$, generate i.i.d.\ random noises $\theta_{i}(k)$ from the uniform distribution $\mathcal{U}(-1,1)$, and randomly select $L$ from the discrete uniform distribution $\mathcal{U}\{1,K_{2}-K_{1}\}$ to satisfy \eqref{eq:subtracted_noises}. The convergence of PR-CPOA is shown in Fig.~\ref*{fig:convergence-PR-CPOA}. The markers on the blue line indicate how many iterations of information dissemination are performed when certain precisions $\epsilon$ are specified. The markers on the orange line represent the actual objective errors $|f_{e}^{*} - f^{*}|$ when those iterations are completed. The relationship between $\log \epsilon$ and $t$ is roughly linear. This phenomenon results from the linear convergence of consensus-type information dissemination. For comparison, we implement DiaDSP\cite{ding2022differentially}, where the step size $\alpha$ and the parameters $d$ and $q$ used in the Laplace distribution are set as $0.05$, $1/3$, and $0.95$, respectively. We run $100$ experiments and plot the average curve. We observe that the solution accuracy of DiaDSP is sacrificed to some extent due to its differentially private mechanism.

The effects of privacy preservation are presented in Fig.~\ref*{fig:dataPrivacy}. This figure demonstrates the relationships between the estimation accuracy $\alpha_{k}$ and the maximum disclosure probability $\beta_{k}$ for a single component $p_{i}^{0}(k)$ when different types of noises $\theta_{i}(k)$ are used. These relationships are explicitly characterized by (\ref{eq:disclosure_prob_single_ele}) in Appendix~\ref{subsec:proof_privacy}. In the experiment, we set $K_{1}=10$, $K_{2}=20$, $p=0.8$, and $\gamma = 10^{-5}$. We consider three types of noises that satisfy uniform, normal, and Laplace distributions. The mean and variance of these noises are $0 \text{ and } 1$, respectively. We observe that $\beta_{k}$ increases with $\alpha_{k}$, which confirms the intuition that a less accurate estimate can be obtained with a higher probability. We also notice that uniformly distributed noises yield the smallest $\beta_{k}$ and thus the most effective preservation of $p_{i}^{0}(k)$. This observation supports the conclusion in \cite{he2018preserving}. Note that the maximum disclosure probability $\beta$ of $p_{i}^{0}$ is the product of all $\beta_{k}$, where $k=1,\ldots,m_{i}+1$ (see (\ref{eq:data_privacy_beta}) in Sec.~\ref{sec:analysis}). The degrees $m_{i}$ of local approximations vary from $20$ to $40$ when the specified precision $\epsilon = 10^{-10}$. Hence, in this case, $\beta$ is an extremely small number given $\alpha_{k}$ and $\beta_{k}$ in the figure.

The robustness of the consensus-type iterations is shown in Fig.~\ref*{fig:robustness}. We consider cases where the time-varying links between agents suffer from different rates of failure, % (i.e., $0\%$, $20\%$ and $50\%$), 
which may result from packet drops or delays. We observe that the iterations still converge, thus ensuring the solution accuracy of the proposed algorithm. Nonetheless, the convergence rates tend to be slower as the link failure rates increase.
% exceeding certain thresholds
    % !TEX root = ..\article.tex
\section{Related Work}\label{sec:literature}
Various distributed optimization algorithms can generally be categorized as primal methods\cite{nedic2009subgradient,shi2015extra,xu2017convergence} and dual-based methods \cite{makhdoumi2017convergence}. The core idea of the primal methods is to combine consensus with gradient-based optimization algorithms, thus achieving consensual convergence in the primal domain.
%There have been extensive researches on designing efficient distributed optimization algorithms, e.g., primal methods\cite{nedic2009subgradient,shi2015extra,xu2017convergence} and dual-based methods \cite{shi2014linear,makhdoumi2017convergence}.
%Thanks to the development of gradient tracking\cite{shi2015extra,di2016next,xu2017convergence,qu2018harnessing}, which enables local agents to approximately track the gradients of the global objective function, the convergence rates of these distributed algorithms can nearly match that of the optimal centralized gradient-based algorithm\cite{boyd2004convex}. %\cite{nesterov2013introductory}. 
%shi2014linear,
The intuition of the dual-based methods is to reformulate the consensus requirement as equality constraints, and then solve the dual problems or carry on primal-dual updates in a distributed manner. % solve the dual problems of the equivalent reformulations
%These carefully constructed dual problems are decoupled, thus easily allowing for the distributed implementations of certain linearly convergent centralized optimization algorithms, e.g., ADMM\cite{shi2014linear,makhdoumi2017convergence}.
For convex problems, distributed algorithms guarantee convergence to globally optimal points; for nonconvex problems, the convergence to first-order stationary solutions is ensured\cite{tatarenko2017non,wai2017decentralized,scutari2019distributed}.

The aforementioned work bridges the gap in convergence behaviors between distributed and centralized optimization algorithms. To deploy these algorithms into applications, some practical issues need to be addressed. These issues include privacy preservation and robustness to allow time-varying, directed communication and asynchronous computations.
% time-varying and directed communication, and asynchronous computations due to lack of coordination, transmission delays, or packet drops.
% mainly centers on bridging the gap in terms of

% issue of privacy preservation
The privacy concern of distributed algorithms has received growing attention. Exchanging exact data can lead to the disclosure of sensitive local objective functions, constraints, and states\cite{han2017differentially}.  
%Conventional approaches are based on the premise that exact local data is exchanged between agents. Nevertheless, if there exist adversaries that intentionally gather certain data necessary for estimation, the sensitive information of objective functions, constraints, and local states can be disclosed \cite{han2017differentially}. To tackle this problem, numerous privacy-preserving consensus and distributed optimization algorithms have been proposed. 
To tackle this problem, a typical approach based on message perturbation is to add random noises to the transmitted data during iterations. The perturbation of the critical data (e.g., states\cite{manitara2013privacy,mo2017privacy,nozari2017differentially,he2018preserving,li2020privacy}, gradients\cite{han2017differentially,hale2018cloud}, step-sizes\cite{lou2018privacy}, and functions\cite{nozari2018differentially,gupta2020preserving}) limits its utility for yielding sensible estimations. %step-sizes\cite{lou2018privacy,ye2020privacy},
%he2019privacy
Some work uses uncorrelated Laplacian or Gaussian noises and develops differentially private consensus\cite{huang2012differentially,nozari2017differentially} and distributed optimization algorithms\cite{han2017differentially,nozari2018differentially,hale2018cloud,zhang2017dynamic,cao2021differentially,ding2022differentially}. The differentially private mechanism offers strong privacy guarantees even against adversaries that own arbitrarily much auxiliary information. Nonetheless, it causes the trade-off between privacy and accuracy\cite{nozari2017differentially,nozari2018differentially}. Other work thus turns to correlated noises and fulfills exact average consensus\cite{manitara2013privacy,mo2017privacy,he2018preserving} and optimization\cite{gupta2020preserving}. % The effects of privacy preservation can be characterized by using the notion of data-privacy\cite{he2018preserving}.
%he2019privacy
There are also methods that utilize state decomposition to achieve complete indistinguishability, provided non-colluding neighbors or private interaction weights exist\cite{wang2019privacy}.
Another typical approach is to apply cryptographic techniques %, and related algorithms include
\cite{ruan2019secure,lu2018privacy}. These methods are suitable if trusted agents or shared keys/secrets exist and the extra computation and communication burdens induced by encryption and decryption are acceptable. %\cite{ruan2019secure,hadjicostis2020privacy,lu2018privacy,zhang2019admm}

% issue of robustness
Another widely investigated issue is robustness.
%Besides the privacy concern, the robustness issues have also been widely investigated.
Time-varying directed communication inhibits the construction of doubly stochastic weight matrices, which are crucial for convergence over undirected graphs. Push-sum-based \cite{nedic2017achieving,scutari2019distributed} and push-pull-based algorithms\cite{pu2020push} overcome this challenge. The former combine the push-sum consensus \cite{kempe2003gossip} with gradient-based methods. %and only require column stochastic weight matrices.
The latter use a row stochastic matrix and a column stochastic matrix to mix solution estimates and gradient trackers, respectively. Algorithms that handle transmission delays include \cite{yang2016distributed}, with the idea of locally fusing the delayed information as soon as it arrives. To achieve asynchronous computations, gossip-type algorithms\cite{xu2017convergence,pu2020push} and those allowing delays and packet drops\cite{tian2020achieving} are developed. %\cite{xin2018linear,pu2020push}
%\cite{wu2017decentralized,tian2020achieving}

Different from the existing work, we exploit polynomial approximation and introduce effective mechanisms to meet practical requirements on privacy and robustness. We show that the proposed algorithm achieves efficient distributed optimization of general nonconvex objectives, and the issues of privacy-accuracy trade-off and step-size selections are avoided.
    % !TEX root = ..\article.tex
\section{Conclusion}\label{sec:conclusion}
We proposed PR-CPOA to solve a class of constrained distributed nonconvex optimization problems, considering privacy preservation and robustness to network imperfections. We achieved exact convergence and effective preservation of the privacy of local objectives by incorporating a new privacy-preserving mechanism for consensus-type iterations. This mechanism utilized the randomness in blockwise insertions of perturbed data, and the privacy degree was explicitly characterized through ($\alpha,\beta$)-data-privacy. We ensured robustness by using the push-sum average consensus protocol as a basis, and we discussed its extensions to maintain the performance when diverse network imperfections exist. We proved that the major benefits brought by the use of polynomial approximation were preserved, and the above demanding requirements were satisfied in the meantime. Future directions include handling noisy evaluations of local objectives and quantized communication.
% for iterations
% diverse imperfections in network communication
% at the same time

    % !TEX root = ..\article.tex
\appendix
\section{Appendix}
\subsection{Proof of Lemma \ref{lem:itr_precision}}\label{subsec:proof_itr_precision}
%\begin{proof}
%    The proof relies on the investigation of consensus iterations and is referred to our online report\cite{he2021private}.
%\end{proof}

%\begin{proof}%
   The proof consists of two steps. First, we prove that the limit value of $p_{i}^{t} \triangleq x_{i}^{t}/y_{i}^{t}(t\in \mathbb{N})$ is $\overbar{p}$, i.e., $\lim_{t\to \infty} p_{i}^{t} = \overbar{p}$. %indeed $\overbar{p}$
   %\begin{equation}\label{eq:consensus_limit}
   %    \lim_{t\to \infty} p_{i}^{t} = \overbar{p}.
   %\end{equation}
   Then, we demonstrate that the meet of the stopping criterion \eqref{eq:itr_stop_rule} is a sufficient condition for \eqref{eq:itr_precision_requirement}. 

   $\bullet$ \textit{Step 1: Proof of the Limit Value}\par
   We consider the $k$-th component of the involved local variables, where $k=1,\ldots,m$. Let
   \begin{alignat*}{2}
       x^{t}&\triangleq [x_{1}^{t}(k),\ldots,x_{N}^{t}(k)]^{T}, &\quad& \theta\triangleq [\theta_{1}(k),\ldots,\theta_{N}(k)]^{T}, \\
       p^{0}&\triangleq [p_{1}^{0}(k),\ldots,p_{N}^{0}(k)]^{T}, && y^{t}\triangleq [y_{1}^{t},\ldots,y_{N}^{t}]^{T}.
   \end{alignat*}
   If the $k$-th components of some $x_{j}^{t}$, $\theta_{j}$, and $p_{j}^{0}(j \in \mathcal{V})$ are null, then they are regarded as $0$ in the expressions. %We investigate the effects of insertions \eqref{eq:insertion_rule} and subtractions \eqref{eq:subtraction} on the accuracy of the consensus-type updates in Algorithm~\ref{alg:alg_overall} as follows.
   
   We first consider the effect of insertions that happened in the first $K_{1}$ iterations. Let $t_{k}$ be the number of iteration when agent $i$ inserts the perturbed state $\tilde{p}_{i}^{0}(k)$. Since $A^{t_{k}}$ is column stochastic, from \eqref{eq:insertion_rule} and \eqref{eq:push_sum_insertion}, we have
   \begin{equation*}
       1^{\top}x^{t_{k}+1} = 1^{\top}A^{t_{k}}x^{t_{k}+} = 1^{\top}x^{t_{k}+} = 1^{\top}x^{t_{k}} + \tilde{p}_{i}^{0}(k).
   \end{equation*}
   %which means that the sum of the components of $x^{t}$ increases by $\tilde{p}_{i}^{0}(k)$. 
   At the end of the $K_{1}$-th iteration, all the agents have inserted their perturbed initial states. Hence,
   \begin{equation*}
       1^{\top}x^{K_{1}} = 1^{\top}x^{0} + \sum_{i\in \mathcal{V}} \tilde{p}_{i}^{0}(k) = \sum_{i\in \mathcal{V}} \tilde{p}_{i}^{0}(k) = 1^{\top}(p^{0} + \theta).
   \end{equation*}
   Then, we focus on the effect of subtractions happened between time $K_{1}+1$ and time $K_{2}$. Suppose that agent $i$ performs its first action of subtractions at the $t_{1}$-th iteration. we use \eqref{eq:push_sum} and the column stochasticity of $A^{t}(t\in \mathbb{N})$ to obtain  %From the column stochasticity of $A^{t}(t\in \mathbb{N})$ and \eqref{eq:push_sum}, it is not difficult to obtain that
   \begin{equation*}
       1^{\top}x^{t_{1}} = 1^{\top}A^{t_{1}-1}x^{t_{1}-1} = 1^{\top}x^{t_{1}-1} = \ldots = 1^{\top}x^{K_{1}}.
   \end{equation*}
   At the $t_{1}$-th iteration, we have
   \begin{equation*}
       1^{\top}x^{t_{1}+1} = 1^{\top}x^{t_{1}} - \delta_{i}(k) = 1^{\top}x^{K_{1}} - \delta_{i}(k).
   \end{equation*}
   %which implies that the sum of the components of $x^{t}$ decreases by $\delta_{i}(k) = \theta_{i}(k)/L$.
   % and that of $y^{t}$ remains the same.
   At the end of the $K_{2}$-th iteration, every agent has completed its $L$ rounds of subtracting the noises. Therefore,
   \begin{equation*}
       1^{\top}x^{K_{2}} = 1^{\top}x^{K_{1}}- 1^{\top}\theta = 1^{\top}p^{0}.
   \end{equation*}
   Since $y_{i}^{t}$ is constantly updated by \eqref{eq:push_sum}, we have
   \begin{equation*}
       1^{\top}y^{K_{2}} = 1^{\top}A^{K_{2}-1} y^{K_{2}-1} = 1^{\top}y^{K_{2}-1} = \ldots = 1^{\top}y^{0}.
   \end{equation*}
   Later on, agents continue to update $x_{i}^{t} \text{ and } y_{i}^{t}$ by \eqref{eq:push_sum}. Based on the convergence of \eqref{eq:push_sum}, we conclude that the exact average can still be achieved, i.e., $\forall k=1,\ldots,m$,
   \begin{equation*}
       \lim_{t\to \infty} p_{i}^{t} = \lim_{t\to \infty} \frac{x_{i}^{t}}{y_{i}^{t}} = \frac{1^{\top}x^{K_{2}}}{1^{\top}y^{K_{2}}} = \frac{1^{\top}p^{0}}{1^{\top}y^{0}} = \frac{\sum_{j=1}^{N} p_{j}^{0}(k)}{N} = \overbar{p}(k).
   \end{equation*}
   %This result holds for any $k=1,\ldots,m$. 
   Therefore, the limit value of $p_{i}^{t}$ is $\overbar{p}$. %, i.e., \eqref{eq:consensus_limit} holds.

   $\bullet$ \textit{Step 2: Proof of the Sufficiency}\par
   Next, we verify the effectiveness of the stopping criterion \eqref{eq:itr_stop_rule}. Note that $p_{i}^{t} = x_{i}^{t}/y_{i}^{t}$, $\forall t\in \mathbb{N}$. The push-sum-consensus-based update of $x_{i}^{t}$ in \eqref{eq:push_sum} can be transformed to
   \begin{equation*}
       p_{i}^{t+1} = \sum_{j=1}^{N} w_{ij}^{t}p_{j}^{t},\quad \text{where } w_{ij}^{t} = \frac{a_{ij}^{t}y_{j}^{t}}{y_{i}^{t+1}}.
   \end{equation*}
   It follows from \eqref{eq:push_sum} and the choice of the weight $a_{ij}^t$ that %$W^{t} \triangleq (w_{ij}^{t})_{N \times N}$ is row stochastic, i.e.,
   $\sum_{j=1}^{N} w_{ij}^{t} = 1$, $w_{ij} \in [0,1]$, $\forall i,j = 1,\ldots,N$.
   %\begin{equation*}
   %    \sum_{j=1}^{N} w_{ij}^{t} = 1, ~~ 0 \leq w_{ij} \leq 1, \quad \forall i,j = 1,\ldots,N,~\forall t.
   %\end{equation*}
   Hence,
   \begin{equation*}
       \begin{split}
           p_{i}^{t+1}(k) &= \sum_{j=1}^{N} w_{ij}^{t} p_{j}^{t}(k) \leq \sum_{j=1}^{N} w_{ij}^{t} \max_{j\in \mathcal{V}} p_{j}^{t}(k) \\
                       &= \max_{j\in \mathcal{V}} p_{j}^{t}(k), \quad \forall k=1,\ldots,m+1,~\forall i\in \mathcal{V}.
       \end{split}
   \end{equation*}
   Let $\displaystyle M^{t}(k) \triangleq \max_{i\in \mathcal{V}} p_{i}^{t}(k)$, $\displaystyle m^{t}(k) \triangleq \min_{i\in \mathcal{V}} p_{i}^{t}(k)$. Then, %It follows that
   \begin{equation*}
       M^{t+1}(k)\leq M^{t}(k), \quad m^{t+1}(k)\geq m^{t}(k).
   \end{equation*}
   It has been proven that $\displaystyle \lim_{t\rightarrow \infty} p_{i}^{t}(k) = \overbar{p}(k)$, $\forall i\in \mathcal{V}$. Hence,
   \begin{equation*}
       \lim_{t\rightarrow \infty} M^{t}(k) = \overbar{p}(k), \quad \lim_{t\rightarrow \infty} m^{t}(k) = \overbar{p}(k).
   \end{equation*}
   Since the sequences of $\big(M^{t}(k)\big)_{t\in \mathbb{N}}$ and $\big(m^{t}(k)\big)_{t\in \mathbb{N}}$ are non-increasing and non-decreasing, respectively, we have
   \begin{equation*}
       m^{t}(k) \leq \overbar{p}(k) \leq M^{t}(k),\quad \forall t\in \mathbb{N}.
   \end{equation*}
   Note that the max/min consensus protocols converge in $U$ iterations. When agents terminate at time $K$, we have
   \begin{equation*}
       r_{i}^{K}(k)-s_{i}^{K}(k) = M^{K'}(k)-m^{K'}(k), 
   \end{equation*}
   where $K'\triangleq K-U$. The meet of \eqref{eq:itr_stop_rule} implies that
   \begin{equation*}
       \begin{split}
       \left|p_{i}^{K}(k)-\overbar{p}(k)\right| &\leq M^{K}(k)-m^{K}(k) \\
       &\leq r_{i}^{K}(k)-s_{i}^{K}(k) \leq \delta,~~\forall i,k.  %\qedhere
       \end{split}
   \end{equation*}
%\end{proof}

\subsection{Proof of Theorem \ref{thm:alg_accuracy_result}}\label{subsec:proof_accuracy}
% \begin{proof}%
   The proof is similar to that of \cite[Theorem~2]{he2020distributed}. We provide a sketch of the main steps here.
   The key idea is to prove the closeness between $p_{i}^{K}(x)$ and $f(x)$ on $X=[a,b]$. Then, their optimal values are also close enough (see \cite[Lemma 1]{he2020distributed}). 
   Note that $p_{i}^{K}(x)$ and $\overbar{p}(x)$ are in the forms of \eqref{eq:cheb_rep} with their coefficients $\{c_{j}\}$ and $\{\overbar{c}_{j}'\}$ stored in $p_{i}^{K}$ and $\overbar{p}$, respectively. It follows from \eqref{eq:itr_precision_requirement} that $\forall x\in [a,b]$,
   \begin{align*}
           \big|&p_{i}^{K}(x) - \overbar{p}(x)\big| = \bigg|\sum_{j=0}^{m} (c_{j}'-\overbar{c}_{j}') T_{j} \Big(\frac{2x-(a+b)}{b-a}\Big) \bigg| \\
           &\leq \sum_{j=0}^{m} \left|c_{j}' - \overbar{c}_{j}' \right| \cdot 1 \leq \sum_{j=0}^{m} \|p_{i}^{K} - \overbar{p}\|_{\infty} \leq \delta(m+1) = \epsilon_{2},
   \end{align*}
   where the first inequality is based on $\left|T_{j}(s)\right|\leq 1, \forall s \in [-1,1]$.
   % Then, we establish the closeness between $\overbar{p}(x)$ and $f(x)$.
   Note that $\overbar{p} $ is the average of all $p_{i}^{0}$. Hence, $\overbar{p}(x)$ is also the average of all $p_{i}(x)$. Based on \eqref{eq:approx_require}, we have
   \begin{equation*}
       \begin{split}
           |&\overbar{p}(x) - f(x)| = \bigg|\frac{1}{N}\sum_{i=1}^{N} \big(p_{i}^{(m_{i})}(x) - f_i(x)\big)\bigg| \\
           &\leq \frac{1}{N} \sum_{i=1}^{N} \big|p_{i}^{(m_{i})}(x) - f_i(x)\big| \leq \frac{1}{N} N \epsilon_{1} = \epsilon_{1}, ~~ \forall x\in [a,b].
       \end{split}
   \end{equation*}
   Given that $\epsilon_{1} = \epsilon_{2} = \epsilon/3$, we have
   \begin{align}\label{eq:approx_func_dist}
   %\begin{equation}\label{eq:approx_func_dist}
   %    \begin{split}
           \big|p_{i}^{K}(x) & - f(x)\big| \leq \left|p_{i}^{K}(x) - \overbar{p}(x)\right| + |\overbar{p}(x) - f(x)| \notag \\
           & \leq \epsilon_{1} + \epsilon_{2} = {\textstyle \frac{2}{3}\epsilon,} \qquad \forall x\in [a,b].
   %    \end{split}
   %\end{equation}
   \end{align}
   Let $p^{*}$ be the optimal value of $p_{i}^{K}(x)$ on $X=[a,b]$. It follows from \cite[Lemma 1]{he2020distributed} that $\left|p^{*}-f^{*}\right| \leq {\textstyle \frac{2}{3}\epsilon}$.
   %\begin{equation*}
   %    \left|p^{*}-f^{*}\right| \leq {\textstyle \frac{2}{3}\epsilon}.
   %\end{equation*}
   Note that $p^{*}\leq f_{e}^{*}\leq p^{*}+\epsilon_{3}= p^{*}+\frac{\epsilon}{3}$. Hence, $f^{*} - {\textstyle \frac{2}{3}\epsilon} \leq p^{*} \leq f_{e}^{*} \leq p^{*} + {\textstyle \frac{\epsilon}{3}} \leq f^{*} + \epsilon$,
   %\begin{equation*}
   %    f^{*} - {\textstyle \frac{2}{3}\epsilon} \leq p^{*} \leq f_{e}^{*} \leq p^{*} + {\textstyle \frac{\epsilon}{3}} \leq f^{*} + \epsilon,
   %\end{equation*}
   which leads to $|f_{e}^{*}-f^{*}|\leq \epsilon$.

   We then characterize the distance between $x_{p}^{*}$ and $x_{f}^{*}$. We consider a small solution accuracy $\epsilon$. %, i.e., a quite strict requirement on the accuracy of the obtained solution.
   It follows from \eqref{eq:approx_func_dist} that
   \begin{equation*}
       f(x_{p}^{*}) \leq p_{i}^{K}(x_{p}^{*}) + {\textstyle \frac{2}{3}\epsilon}, \quad p_{i}^{K}(x_{p}^{*}) \leq p_{i}^{K}(x_{f}^{*}) \leq f(x_{f}^{*}) + {\textstyle \frac{2}{3}\epsilon},
   \end{equation*}
   which implies that $f(x_{p}^{*}) \leq f(x_{f}^{*}) + \frac{4}{3}\epsilon$.
   %\begin{equation*}
   %    f(x_{p}^{*}) \leq f(x_{f}^{*}) + {\textstyle \frac{4}{3}\epsilon}.
   %\end{equation*}
   Hence, $x_{p}^{*}$ falls in the sublevel set $S = \{x\in X| f(x) \leq f(x_{f}^{*}) + \frac{4}{3}\epsilon\}$ of $f(x)$.
   %\begin{equation*}
   %    S = \left\{x\in X| f(x) \leq f(x_{f}^{*}) + {\textstyle \frac{4}{3}\epsilon}\right\}.
   %\end{equation*}
   Therefore, we have $|x_{p}^{*} - x_{f}^{*}| \leq \diam(S)$.
   % \begin{equation*}
   %     |x_{p}^{*} - x_{f}^{*}| \leq \diam(S). \qedhere
   % \end{equation*}
% \end{proof}

\subsection{Proof of Theorem \ref{thm:privacy}}\label{subsec:proof_privacy}
%\begin{proof}
    We first consider the estimation of $p_{i}^{0}(k)$, where $k=1,\ldots,m_{i}+1$. Suppose that at the $t_{k}$-th iteration, agent $i$ inserts the perturbed state $\tilde{p}_{i}^{0}(k)$ by \eqref{eq:insertion_rule}. The estimation $\hat{p}_{i}(k)$ of $p_{i}^{0}(k)$ can be calculated at three types of time, i.e., before $t_{k}$, at $t_{k}$, and after $t_{k}$. We discuss these cases as follows.

    % \noindent $\bullet$
    $\bullet$ \textit{Case 1:} At time $t<t_{k}$, $\tilde{p}_{i}^{0}(k)$ has not been inserted yet. What the adversaries have collected are either null values or combinations of the perturbed states of agent $i$'s neighbors. Since there is not any available information on $p_{i}^{0}(k)$ that serves as a basis for estimation, by \eqref{eq:est_acc_no_info}, we have
    \begin{equation*}
        \PR{|\hat{p}_{i}(k)-p_{i}^{0}(k)|\leq \alpha_{k} \big|\mathcal{I}_{i}^{t}} \leq \gamma.
    \end{equation*}
    $\bullet$ \textit{Case 2:} At time $t=t_{k}$, $\tilde{p}_{i}^{0}(k)$ is inserted. By Assumption \ref{assump:bound_ability_adversary}, the probability that the adversaries acquire the full knowledge of $\mathit{I}_{i}^{\textrm{in},t_{k}-1}$ is not more than $p$. If this is the case, based on \eqref{eq:insertion_rule} and \eqref{eq:push_sum_insertion}, they can easily calculate $\tilde{p}_{i}^{0}(k)$ by
    \begin{equation}\label{eq:recover_state}
        \tilde{p}_{i}^{0}(k) = x_{i}^{t_{k}+}(k) - \sum_{j \in \mathcal{N}_{i}^{\text{in},t_{k}-1}} a_{ij}^{t_{k}-1} x_{j}^{(t_{k}-1)+}(k).
    \end{equation}
    Note that $\tilde{p}_{i}^{0}(k) = {p}_{i}^{0}(k) + \theta_{i}(k)$.
    %\begin{equation*}
    %    \tilde{p}_{i}^{0}(k) = {p}_{i}^{0}(k) + \theta_{i}(k).
    %\end{equation*}
    Hence, after an estimation $\hat{\theta}_{i}(k)$ of $\theta_{i}(k)$ is obtained, $\hat{p}_{i}(k)$ is calculated by $\hat{p}_{i}(k) = \tilde{p}_{i}^{0}(k) - \hat{\theta}_{i}(k)$.
    %\begin{equation*}
    %    \hat{p}_{i}(k) = \tilde{p}_{i}^{0}(k) - \hat{\theta}_{i}(k).
    %\end{equation*}
    Therefore, for any estimation $\hat{p}_i(k)$, we have
    \begin{align}\label{eq:privacy_tk_1}
        \Pr &\big\{|\hat{p}_{i}(k)-p_{i}^{0}(k)|\leq \alpha_{k} \big|\mathcal{I}_{i}^{t_{k}}\big\} \notag \\
        & = \PR{|\hat{\theta}_{i}(k)-\theta_{i}(k)|\leq \alpha_{k} \big|\mathcal{I}_{i}^{t_{k}}} \notag \\
        & = \PR{\theta_{i}(k)\in [\hat{\theta}_{i}(k)-\alpha_{k},\hat{\theta}_{i}(k)+\alpha_{k}] \big|\mathcal{I}_{i}^{t_{k}}} \notag \\
        & = \int_{\hat{\theta}_{i}(k)-\alpha_{k}}^{\hat{\theta}_{i}(k)+\alpha_{k}} g_{\theta_{i}(k)}(y) \ud y \notag \\
        & \leq \max_{\nu \in \Theta} \int_{\nu-\alpha_{k}}^{\nu+\alpha_{k}} g_{\theta_{i}(k)}(y) \ud y,
    \end{align}
    where $\hat{\theta}_{i}(k) \in \Theta$.
    However, if the adversaries can only access part of $\mathit{I}_{i}^{t_{k}-1}$, they are unable to calculate $x_{i}^{t_{k}}(k)$ by \eqref{eq:push_sum_insertion} and then recover $\tilde{p}_{i}^{0}(k)$ by \eqref{eq:recover_state}. Note that
    \begin{equation*}
        x_{i}^{t_{k}+}(k) = x_{i}^{t_{k}}(k) + \tilde{p}_{i}^{0}(k) = x_{i}^{t_{k}}(k) + \theta_{i}(k) + {p}_{i}^{0}(k).
    \end{equation*}  
    Hence, in this case, they need to obtain an estimation $\hat{\eta}_{i}(k)$ of $x_{i}^{t_{k}}(k) + \theta_{i}(k)$, and then calculate $\hat{p}_{i}(k)$ by $\hat{p}_{i}(k) = x_{i}^{t_{k}+}(k) - \hat{\eta}_{i}(k)$.
    %\begin{equation*}
    %    \hat{p}_{i}(k) = x_{i}^{t_{k}+}(k) - \hat{\eta}_{i}(k).
    %\end{equation*}
    According to \eqref{eq:push_sum_insertion}, $x_{i}^{t_{k}}(k)$ is a linear combination of the states $x_{j}^{(t_{k}-1)+}$ for $j\in \mathcal{N}_{i}^{\text{in},t_{k}-1}$. These states depend on some $\tilde{p}_{l}^{0}(k)$ and thus also some $\theta_{l}(k)$, where $l\in \mathcal{V}$. Note that the adversaries only have partial knowledge of $\mathit{I}_{i}^{\textrm{in},t_{k}-1}$ and know part of these states. Hence, there exist certain independent random variables, i.e., $\theta_{l}(k)$, of which the adversaries do not own any prior or relevant knowledge. As a result, by \eqref{eq:est_acc_no_info}, it is hard to estimate $x_{i}^{t_{k}}(k)$ with high precision. It follows that
    \begin{align}\label{eq:privacy_tk_2}
        & \Pr \big\{|\hat{p}_{i}(k)-p_{i}^{0}(k)|\leq \alpha_{k}\big|\mathcal{I}_{i}^{t_{k}} \big\} \notag \\
        & ~= \PR{\big|\hat{\eta}_{i}(k)-(x_{i}^{t_{k}}(k)+\theta_{i}(k))\big|\leq \alpha_{k}\big|\mathcal{I}_{i}^{t_{k}}} \notag \\
        & ~\leq \PR{\hat{\eta}_{i}(k)-x_{i}^{t_{k}}(k) \in [\theta_{i}(k)-\alpha_{k}, \theta_{i}(k)+\alpha_{k}] \big|\mathcal{I}_{i}^{t_{k}},\theta_{i}(k)} \notag \\
        & ~\leq \gamma,
    \end{align}
    % We combine \eqref{eq:privacy_tk_1} and \eqref{eq:privacy_tk_2} and know that for any estimation $\hat{p}_i(k)$ of $p_i^0(k)$
    Combining \eqref{eq:privacy_tk_1} and \eqref{eq:privacy_tk_2}, for any estimation $\hat{p}_i(k)$ of $p_i^0(k)$, we have
    \begin{equation}\label{eq:privacy_tk_result}
    \begin{split}
        \max_{\hat{p}_i(k)}~ & \Pr \big\{|\hat{p}_{i}(k)-p_{i}^{0}(k)|\leq \alpha_{k}\big|\mathcal{I}_{i}^{t_{k}} \big\} \\
            &\leq p\max_{\nu \in \Theta} \int_{\nu-\alpha_{k}}^{\nu+\alpha_{k}} g_{\theta_{i}(k)}(y) \ud y + \gamma \triangleq h_{i}(\alpha_{k}).
    \end{split}
    \end{equation}
    $\bullet$ \textit{Case 3:} At time $t > t_{k}$, the adversaries can estimate $p_{i}^{0}(k)$ either by the same rule that is adopted at time $t = t_{k}$ or by the new rule based on the new information. In the former case, we still obtain \eqref{eq:privacy_tk_result}. We now discuss the latter case in detail. We first consider the time $t = t_{k}+1$. Note that
    \begin{align}\label{eq:update_rule_beyond}
        & \frac{x_{i}^{(t_{k}+1)+}(k)}{a_{ii}^{t_{k}}} = \frac{x_{i}^{t_{k}+1}(k)}{a_{ii}^{t_{k}}} \notag \\
            &\phantom{==} = x_{i}^{t_{k}+}(k) + \frac{1}{a_{ii}^{t_{k}}}\Big(\sum_{j \in \mathcal{N}_{i}^{\text{in},t_{k}}\setminus \{i\}} a_{ij}^{t_{k}} x_{j}^{t_{k}+}(k) - \tau_{i,t_{k}+1}(k)\Big) \notag \\
            &\phantom{==} = p_{i}^{0}(k) + \theta_{i}(k) + x_{i}^{t_{k}}(k) \notag \\
            &\phantom{===} + \frac{1}{a_{ii}^{t_{k}}}\Big(\sum_{j \in \mathcal{N}_{i}^{\text{in},t_{k}}\setminus \{i\}} a_{ij}^{t_{k}} x_{j}^{t_{k}+}(k) - \tau_{i,t_{k}+1}(k)\Big) \notag \\
            &\phantom{==} = p_{i}^{0}(k) + \theta_{i}(k) + \theta_{i}'(k),
    \end{align}
    where $\tau_{i,t}(k) = \zeta_{i}(k)$ if noises are subtracted at time $t$, and $\tau_{i,t}(k) = 0$ otherwise. If the full knowledge of $\mathit{I}_{i}^{\textrm{in},t_{k}}$ is available, the adversaries can not only collect all the $x_{j}^{t_{k}+}$ for $j\in \mathcal{N}_{i}^{\text{in},t}$, but also accurately infer $\tau_{i,t_{k}+1}(k)$ by
    \begin{equation*}
        \tau_{i,t_{k}+1}(k) = \sum_{j \in \mathcal{N}_{i}^{\text{in},t_{k}}} a_{ij}^{t_{k}} x_{j}^{t_{k}+}(k) - x_{i}^{(t_{k}+1)+}(k). 
    \end{equation*}
    Hence, $\theta_{i}'(k)$ is a deterministic constant. In this case, by using \eqref{eq:update_rule_beyond}, we still have
    \begin{align*}
        \Pr &\big\{|\hat{p}_{i}(k)-p_{i}^{0}(k)|\leq \alpha_{k} \big|\mathcal{I}_{i}^{t_{k}+1}\big\} \\
        & = \PR{|\hat{\theta}_{i}(k)-\theta_{i}(k)|\leq \alpha_{k} \big|\mathcal{I}_{i}^{t_{k}+1}}.
        % & = \int_{\hat{\theta}_{i}(k)-\alpha_{k}}^{\hat{\theta}_{i}(k)+\alpha_{k}} g_{\theta_{i}(k)}(y) \ud y \notag \\
        % & \leq \max_{\nu \in \Theta} \int_{\nu-\alpha_{k}}^{\nu+\alpha_{k}} g_{\theta_{i}(k)}(y) \ud y.
    \end{align*}
    Next, we analyze the disclosure probability of $\theta_{i}(k)$ given $\mathcal{I}_{i}^{t_{k}+1}$. The newly available information, i.e., the subtracted noise $\zeta_{i}(k)$, allows for another means of inferring $\theta_{i}(k)$. 
    % The second means is to infer $\theta_{i}(k)$ based on the available subtracted noise $\zeta_{i}(k)$.
    We now show that the resulting disclosure probability is rather small when $L_i$ is drawn from an unknown distribution. Note that $\zeta_{i}(k) = \theta_{i}(k)/L_i > \alpha_{k}$. Hence,
    \begin{align}
        \Pr &\big\{|\hat{p}_{i}(k)-p_{i}^{0}(k)|\leq \alpha_{k} \big|\mathcal{I}_{i}^{t_{k}+1}\big\} \notag \\
        & = \PR{|\hat{\theta}_{i}(k)-\theta_{i}(k)|\leq \alpha_{k} \big|\zeta_{i}(k)} \notag \\
        & = \PR{|\hat{L}_i - L_i|\cdot\zeta_{i}(k)\leq \alpha_{k}\big|\zeta_{i}(k)} \notag \\
        & = \PR{\hat{L}_i = L_i\big|\zeta_{i}(k)} \leq \gamma,
    \end{align}
    where $\hat{L}_i$ is any estimation of $L_i$, and the last inequality follows from \eqref{eq:est_acc_no_info}. Thus, the disclosure probability will not exceed the upper bound in \eqref{eq:privacy_tk_1}, i.e.,
    \begin{equation}\label{eq:privacy_tk+1_1}
        \begin{split}
            \Pr &\big\{|\hat{p}_{i}(k)-p_{i}^{0}(k)|\leq \alpha_{k} \big|\mathcal{I}_{i}^{t_{k}+1}\big\} \\
            & = \PR{|\hat{\theta}_{i}(k)-\theta_{i}(k)|\leq \alpha_{k} \big|\mathcal{I}_{i}^{t_{k}+1}} \\
            & \leq \max_{\nu \in \Theta} \int_{\nu-\alpha_{k}}^{\nu+\alpha_{k}} g_{\theta_{i}(k)}(y) \ud y.
        \end{split}
    \end{equation}
    If the full knowledge of $\mathit{I}_{i}^{\textrm{in},t}$ is unavailable, then $\theta_{i}'(k)$ contains those independent random variables whose relevant information is unknown to the adversaries. Specifically, if $t_{k}+1 \leq K_{1}$, then those variables refer to certain added noises $\theta_{l}(k)$ that are included in $x_{l}^{t_{k}+}(k)$, where $l\in \mathcal{V}$. Otherwise, those variables refer to certain subtracted noises $\zeta_{l}(k)$ for some $l\in \mathcal{V}$. Thus, it follows from \eqref{eq:est_acc_no_info} that
    \begin{equation}\label{eq:privacy_tk+1_2}
        \PR{|\hat{p}_{i}(k)-p_{i}^{0}(k)|\leq \alpha_{k} \big|\mathcal{I}_{i}^{t_{k}+1}}\leq \gamma.
    \end{equation}
    Combining \eqref{eq:privacy_tk+1_1} and \eqref{eq:privacy_tk+1_2}, for any estimation $\hat{p}_i(k)$ of $p_i^0(k)$, % we have
    \begin{equation}\label{eq:privacy_tk+1_result}
    \begin{split}
        \max_{\hat{p}_i(k)}~ & \Pr \big\{|\hat{p}_{i}(k)-p_{i}^{0}(k)|\leq \alpha_{k}\big|\mathcal{I}_{i}^{t_{k}+1} \big\} \\
            &\leq p\max_{\nu \in \Theta} \int_{\nu-\alpha_{k}}^{\nu+\alpha_{k}} g_{\theta_{i}(k)}(y) \ud y + \gamma = h_{i}(\alpha_{k}).
    \end{split}
    \end{equation}
    A similar analysis can be performed for other arbitrary $t \geq t_{k}+1, t\in \mathbb{N}$. 
    However, for $t\geq K_{2}$, there exists an extreme case where the adversaries successfully obtain the full knowledge of $I_{i}^{\text{in},t}$ at time $t=t_{k}-1$ and also from time $t=K_{1}+1$ to time $t=K_{2}$. In this case, they can not only calculate $\tilde{p}_{i}^{0}(k)$ by \eqref{eq:recover_state}, but also acquire $\tau_{i,t}(k)$ and perfectly infer $\theta_{i}(k)$ by $\theta_{i}(k) = \sum_{t=K_{1}+1}^{K_{2}} \tau_{i,t}(k)$.
    %\begin{equation*}
    %    \theta_{i}(k) = \sum_{t=K_{1}+1}^{K_{2}} \tau_{i,t}(k).
    %\end{equation*}
    Hence, the exact value of $p_{i}^{0}(k)$ can be obtained, and
    \begin{align*}
        \Pr &\big\{|\hat{p}_{i}(k)-p_{i}^{0}(k)|\leq \alpha_{k} \big|\mathcal{I}_{i}^{K_{2}}\big\} \\
        & = \PR{|\hat{\theta}_{i}(k)-\theta_{i}(k)|\leq \alpha_{k} \big|\mathcal{I}_{i}^{K_{2}}} = 1.
    \end{align*}
    The probability that such an extreme case happens is not more than $p^{K_{2}-K_{1}+1}$. Thus, for any $k=1,\ldots,m_{i}+1\text{ and } t\in \mathbb{N}$, % we have
    \begin{equation}\label{eq:disclosure_prob_values}
            \max_{\hat{p}_i(k)}~ \Pr \big\{|\hat{p}_{i}(k)-p_{i}^{0}(k)|\leq \alpha_{k}\big|\mathcal{I}_{i}^{t} \big\} \leq \beta_{k},
    \end{equation}
    where $\beta_k$ is given by \eqref{eq:disclosure_prob_single_ele}.
    % \begin{align}\label{eq:disclosure_prob_single_ele}
    %     \beta_{k} = \big(1-p^{K_{2}-K_{1}+1}\big) h_{i}(\alpha_{k}) + p^{K_{2}-K_{1}+1}.
    % \end{align}
    Since $h_{i}(\alpha_{k})\leq p+\gamma < 1$, $\beta_{k}$ is larger than the RHS of \eqref{eq:privacy_tk_result}. 

    Finally, we consider the inference on whether $p_{i}^{0}(k)$ is null for $k=m_{i}+2,\ldots,m+1$, i.e., whether $p_{i}^{0}$ is an $(m_{i}+1)$-dimensional vector. Note that there is no action of insertions or subtractions corresponding to the aforementioned components. Hence, the adversaries will not find any inconsistency between $x_{i}^{(t+1)+}(k)$ and $\sum_{j \in \mathcal{N}_{i}^{\text{in},t}} a_{ij}^{t} x_{j}^{t+}(k)$, where $t\in \mathbb{S}_{t}$. Let this event be denoted by $A$. Once it occurs, the adversaries need to decide between the following two hypotheses, i.e., $\mathcal{H}_{0}: p_{i}^{0}(k)$ is null, and $\mathcal{H}_{1}: p_{i}^{0}(k)$ is a nonzero number.
    %\begin{equation*}
    %    \mathcal{H}_{0}: p_{i}^{0}(k) \text{ is null,} \qquad \mathcal{H}_{1}: p_{i}^{0}(k) \text{ is a nonzero number.}
    %\end{equation*}
    Based on the algorithmic design, we have $\PR{A|\mathcal{H}_{0}} = 1, \PR{A|\mathcal{H}_{1}} = (1-p)^{L_i+1}$.
    %\begin{equation*}
    %    \PR{A|\mathcal{H}_{0}} = 1, \qquad \PR{A|\mathcal{H}_{1}} = (1-p)^{L+1}.
    %\end{equation*}
    It follows from the maximum likelihood rule that the adversaries will always choose $\mathcal{H}_{0}$ when $A$ occurs. The probability that they successfully decide that $p_{i}^{0}(k)$ is null for $k=m_{i}+2,\ldots,m+1$ equals
    \begin{equation}\label{eq:disclosure_prob_null}
        \PR{m_{i}+1 \leq k-1} = F_{m_{i}|\mathcal{I}_{i}^{t}}(k-2).
    \end{equation}
    Combining \eqref{eq:disclosure_prob_values} and \eqref{eq:disclosure_prob_null}, we arrive at \eqref{eq:data_privacy_beta}.
    %we have
    %\begin{align*}%\label{eq:privacy_result}
    %    \max_{\hat{p}_i}~ & \Pr \big\{\|\hat{p}_{i} - p_{i}^{0}\|_{1}\leq \alpha|\mathcal{I}_{i}^{t}\big\} \notag \\ 
    %        & = \prod_{k=1}^{m_{i}+1} \max_{\hat{p}_i(k)~} \PR{|\hat{p}_{i}(k)-p_{i}^{0}(k)|\leq \alpha_{k}\big|\mathcal{I}_{i}^{t}} \notag \\
    %        & \qquad \cdot \prod_{k=m_{i}+2}^{m+1} \PR{m_{i}+1 \leq k-1} \leq 
    %        \beta  %\notag \\
    %        %& \leq \prod_{k=1}^{m_{i}+1} \beta_{k} \cdot \prod_{k=m_{i}+2}^{m} F_{m_{i}|\mathcal{I}_{i}^{t}}(k-2) = \beta. %\qedhere
    %\end{align*}
    % where $\beta$ is given by (\ref{eq:data_privacy_beta}).
%\end{proof}

\subsection{Proof of Theorem~\ref{thm:complexity}}\label{subsec:proof_complexity}
%\begin{proof}
    Note that the evaluations of local objective functions (i.e., queries of the zeroth-order oracle) are only performed in the stage of initialization, and the primal-dual interior-point method\cite{boyd2004convex} is used to solve the reformulated SDP in the stage of polynomial optimization. By referring to the proof of \cite[Theorem 6]{he2020distributed}, we know that for every agent, the orders of evaluations of local objective functions and primal-dual iterations are of $\bigO{m}$ and $\bigO{\sqrt{m}\log \frac{1}{\epsilon}}$, respectively. %, where $m$ is the maximum degree of local approximations. 
    Also, the orders of the required flops of these two stages are of $\bigO{m\cdot\max(m,F_{0})}$ and $\bigO{m^{4.5} \log \frac{1}{\epsilon}}$, respectively.

    In the stage of information dissemination, the blockwise insertions of vectors and the subtractions of noises are completed in finite time, i.e., in $K_{2}$ iterations. Since the consensus-type protocol converges geometrically, the order of the total number of iterations (i.e., inter-agent communication) is of
    \begin{equation*}    
        K_{2} + \mathcal{O}\Big(\log \frac{1}{\delta}\Big) = \mathcal{O}\Big(\log \frac{1}{\delta}\Big) = \bigO{\log \frac{m}{\epsilon}},
    \end{equation*}
    where the required precision $\delta$ is given by \eqref{eq:itr_stop_rule}. The order of flops needed in this stage is of $\bigO{m \log \frac{m}{\epsilon}}$. The results in the theorem follow from the above analysis.
%\end{proof}

    \balance
    \bibliographystyle{IEEEtran}
    \bibliography{article}

\end{document}